\documentclass{article}

\usepackage{amsmath,amssymb,stmaryrd,color,bm,nicefrac,amsthm}
\usepackage{color}    
\usepackage{here}
\usepackage{graphicx, subfigure,stmaryrd,enumerate}
\usepackage{proof}

\usepackage[
    colorlinks=true, 
    citecolor=red,
    linkcolor=orange,
    urlcolor=blue]{hyperref}

\usepackage{tikz}

\newcommand{\dom}{\mathop{\rm dom}}
\newcommand{\Regular}[1]{\Omega_{#1}}
\newcommand{\mlang}{\mathcal L_\Box}

\newcommand{\zfc}{{\rm ZFC}}
\newcommand{\Ord}{{\sf Ord}}

\newcommand{\lb}{\llbracket}
\newcommand{\rb}{\rrbracket}
\newcommand{\gl}{{\sf GL}}
\newcommand{\glp}{{\sf GLP}}

\newcommand{\cf}[1]{\mathop{\text{\it cf}} #1}

\newcommand{\jp}[1]{}

\newcommand{\ic}[2]{{#2}_{+#1}}
\newcommand{\ico}[1]{\mathcal{I}_{#1}}

\newtheorem{theorem}{Theorem}[section]

\newtheorem{exm}[theorem]{Example}

\newtheorem{question}{Question}
\newtheorem{defi}[theorem]{Definition}

\newtheorem{lemma}[theorem]{Lemma}
\newtheorem{prop}[theorem]{Proposition}
\newtheorem{cor}[theorem]{Corollary}

\usepackage{authblk}

\title{Strong Completeness of Provability Logic for Ordinal Spaces}
\author{Juan P. Aguilera}
\affil{Vienna University of Technology, Austria.}
\author{David Fern\'{a}ndez-Duque}
\affil{University of Toulouse, France, and Instituto Tecnol\'ogico Aut\'onomo de M\'exico, Mexico.}
\date{}                     
\begin{document}
\maketitle
\begin{abstract}
Abashidze and Blass independently proved that the modal logic $\gl$ is complete for its topological interpretation over any ordinal greater than or equal to $\omega^\omega$ equipped with the interval topology. Icard later introduced a family of topologies $\mathcal I_\lambda$ for $\lambda < \omega$, with the purpose of providing semantics for Japaridze's polymodal logic $\glp_\omega$. Icard's construction was later extended by Joosten and the second author to arbitrary ordinals $\lambda \geq \omega$. 

We further generalize Icard topologies in this article. Given a scattered space $\mathfrak X = (X, \tau)$ and an ordinal $\lambda$, we define a topology $\ic \lambda \tau$ in such a way that $\ic 0 \tau$ is the original topology $\tau$ and $\ic \lambda \tau$ coincides with $\mathcal I_\lambda$ when $\mathfrak X$ is an ordinal endowed with the left topology.

We then prove that, given any scattered space $\mathfrak X$ and any ordinal $\lambda>0$ such that the rank of $(X, \tau)$ is large enough, $\gl$ is strongly complete for $\ic \lambda \tau$. One obtains the original Abashidze-Blass theorem as a consequence of the special case where $\mathfrak X=\omega^\omega$ and $\lambda=1$. 
\end{abstract}

{\maketitle}


\section{Introduction}
The study of formalized provability can be traced back to the seminal 1931 paper by Kurt G\"odel \cite{godel31}, where he first showed how any mathematical theory that is able to describe the arithmetic of natural numbers, such as Peano Arithmetic (PA), also possesses the expressive power to reason about provability within said theory. For this, he first noted that sentences \emph{about} numbers can be encoded using numbers themselves. G\"odel then showed that we can define in the language of arithmetic a formula $B(n)$ that holds true if, and only if, $n$ encodes a formula that is provable in the theory. 

He also noted that many statements about provability may be captured in the language $\mlang$ of modal logic. Formulae in $\mlang$ are built from a countable set of propositional variables $\mathbb{P}$ by applying Boolean connectives and the operator $\Box$, together with its dual $\Diamond$. The provability logic $\gl$ has as axioms the schemata
\[\begin{array}{lllll}
\text{\bf K:}&\Box(\varphi \rightarrow \psi) \rightarrow (\Box\varphi \rightarrow \Box \psi) &\hfill&
{\text{\bf L\"ob:}}& \Box(\Box \varphi \rightarrow \varphi) \rightarrow \Box \varphi
\end{array}\]
and the inference rules \emph{modus ponens} and \emph{necessitation}: $\frac {\varphi}{\Box \varphi}$.

The intended interpretation for provability logic is to read $\Box \varphi$ as ``$\varphi$ is provable in $T$'' for some formal theory $T$. However, other interpretations for $\gl$ are possible. As was proved by K. Segerberg \cite{segerberg1971}, $\gl$ is complete with respect to its Kripke semantics, albeit not \emph{strongly} complete. L. Esakia \cite{esakia1981} showed that $\gl$ is also complete with respect to its topological interpretation, and indeed it is strongly complete. M. Abashidze \cite{abashidze1985} and A. Blass \cite{blass1990} improved Esakia's theorem by showing that $\gl$ is also complete with respect to a single space; namely, an ordinal number with what is known as the \emph{order topology}. It follows from the results of L. Beklemishev and D. Gabelaia \cite{topocompletenessofglp} and the second author \cite{polytopologies} that $\gl$ is also complete for a large collection of spaces introduced by T. Icard \cite{icardglp} and some of its generalizations \cite{modelsofglp}. These topologies find natural interpretations in the construction of models of the polymodal logics $\glp_\Lambda$ with $\Lambda$ modalities, although technical difficulties have stopped these results from being extended to uncountable $\Lambda$. This work partially solves these issues. 

Our goal is to show how one can start with any topological space $\mathfrak X=(X,\tau)$ that validates the axioms of $\gl$ and modify it slightly so as to obtain a class of spaces of the form $\ic\lambda{\mathfrak X}=(X,\ic\lambda\tau)$ with $\lambda$ an ordinal, such that $\gl$ is strongly complete with respect to $\ic\lambda{\mathfrak X}$ for each $\lambda<\Lambda$, where $\Lambda$ depends on $\mathfrak X$ (although it may well be arbitrarily large). We call the spaces $\ic\lambda{\mathfrak X}$ thus obtained {\em generalized Icard spaces.} We expect the present work to be of interest, as strong completeness of provability logic for single spaces is novel, and generalized Icard spaces give a wide class of spaces which provide natural and easily described topological models for $\gl$.

\subsection*{Outline} The article consists of six main sections. In Section \ref{SecOrd}, we review some facts about ordinal arithmetic and transfinite iterations of normal and initial functions. Section \ref{SecTopo} reviews the topological semantics of modal logic. In Section \ref{SecRel}, we introduce a class of models for provability logic called \emph{$\omega$-bouquets} and prove the strong completeness of $\gl$ with respect to this class. In Section \ref{SecIc}, we introduce \emph{generalized Icard spaces} and our main result, strong completeness of $\gl$, is proved in Section \ref{SecGL}. We then finish with some concluding remarks.

\section{Ordinal numbers}\label{SecOrd}

We will be dealing extensively with ordinal numbers and so we assume knowledge of elementary ordinal arithmetic, although we will review some basic facts, as well as some operations that will be of great importance. For a thorough review, we refer the reader to, for example, \cite{jech}.

\begin{lemma}\
\label{ordinal arithmetic}
\begin{enumerate}[(i)]
\item If $\alpha < \beta$, there exists a unique ordinal $\gamma$ such that $\alpha + \gamma = \beta$. We denote this $\gamma$ by $-\alpha + \beta$.
\item\label{OrdArII} For all $\gamma$, if $\alpha > 0$, there exist unique $\beta$ and unique $\rho < \alpha$ such that $\gamma < \alpha\cdot\beta + \rho$.
\item For all nonzero $\xi$, there exist ordinals $\alpha$ and $\beta$ such that $\xi = \alpha + \omega^\beta$. Such a $\beta$ is unique. We denote it by $\ell\xi$ and call it the \emph{end logarithm} of $\xi$.
\end{enumerate}
\end{lemma}

By a $\lambda$-\emph{sequence} of ordinals, we mean an ordinal-valued function with domain $\lambda$. As usual, if $\lambda$ is a limit ordinal and $(\xi_\eta)_{\eta < \lambda}$ is a $\lambda$-sequence, we define $\xi=\lim_{\eta\to\lambda}\xi_\eta$ if for every $\zeta<\xi$ there is $\delta<\lambda$ such that $\xi_\eta\in (\zeta,\xi]$ whenever $\eta>\delta$. In particular, if $(\xi_\eta)_{\eta < \lambda}$ is non-decreasing, then $\lim_{\eta \to \lambda} \xi_\eta = \bigcup\{\xi_\eta \colon \eta \in \lambda\}.$ A (class) function $f$ on the ordinals is \emph{continuous} if $f(\lambda) = \lim_{\eta \to \lambda} f(\eta)$ for each limit $\lambda$. A function that is both continuous and increasing is called \emph{normal}. An important example is given as follows:

\begin{defi}[Exponential function]
The \emph{exponential function} is the normal function $e\colon \Ord \to \Ord$ given by \mbox{$\xi \mapsto -1 + \omega^\xi$}.
\end{defi}

When $f\colon X\to X$ is a function, it is natural and often useful to ask whether $f$ has {\em fixed points,} i.e., solutions to the equation $x=f(x)$. In particular, normal functions have many fixed points:

\begin{prop}\label{fixedpoints}
Every normal function on $\Ord$ has arbitrarily large fixed points.
\end{prop}

The first ordinal $\alpha$ such that $\alpha = \omega^\alpha$ is the limit of the $\omega$-sequence\linebreak $(\omega, \omega^\omega, \omega^{\omega^\omega}, ...)$, and is usually called $\varepsilon_0$. In general, we call an \emph{epsilon number} any nonzero fixed point of the exponential function, and \emph{epsilon function} the function given by $\alpha \mapsto \varepsilon_\alpha$, assumed to be increasing and to have as range all epsilon numbers. It is easily proved that the epsilon function is normal and $\varepsilon_{\alpha+1} = \sup\{\omega^{\varepsilon_{\alpha}+1}, \omega^{\omega^{\varepsilon_{\alpha}+1}}, ...\}$ for each $\alpha$.

Recall that a non-zero ordinal number $\xi$ is said to be \emph{additively indecomposable} if $\beta+\gamma < \xi$ for all $\beta, \gamma < \xi$. Also, $\xi$ is said to be \emph{multiplicatively indecomposable} if $\beta\gamma < \xi$ for all $\beta, \gamma < \xi$. The following facts are well-known:

\begin{lemma}\
\label{indecomposable}
\begin{enumerate}
\item\label{indecomposableOne} An ordinal is additively indecomposable if, and only if, it is of the form $\omega^\rho$.
\item An ordinal is multiplicatively indecomposable if, and only if, it is of the form $\omega^{\omega^\rho}$.
\end{enumerate}
\end{lemma}

In \cite{hyperations}, the second author and J. Joosten introduced families of transfinite iterations of functions called \emph{hyperations} and \emph{cohyperations}. Although we will not provide a detailed discussion, we consider two examples which have remarkable applications to provability logic, as will be shown later on.

\begin{defi}[Hyperexponential functions]
The \emph{hyperexponentials} $(e^\zeta)_{\zeta \in \Ord}$ are the unique family of normal functions that satisfy
\begin{enumerate}[(i)]
\item\label{DefHypI} $e^1 = e$,
\item\label{DefHypII} $e^{\alpha + \beta} = e^\alpha \circ e^\beta$ for all $\alpha$ and $\beta$, and
\item if $(f^\xi)_{\xi\in\Ord}$ is a family of functions satisfying \ref{DefHypI} and \ref{DefHypII}, then for all $\alpha,\beta\in\Ord$, $e^\alpha\beta\leq f^\alpha\beta$.
\end{enumerate}
\end{defi}

This definition is proved to be correct in \cite{hyperations}. From now on, to ease notation, we will sometimes omit the symbol `$\circ$' for composition of functions, as well as parentheses. In order to compute hyperexponentials, we use the following lemma:

\begin{lemma}[Recursive hyperexponentiation]\
\label{recexp}
\begin{enumerate}[(i)]
\item  $e^0$ is the identity;
\item  $e^\xi 0 = 0$ for all $\xi$;
\item  \label{RecHypIII}for any $\xi$ and any limit $\lambda$, $e^\xi\lambda = \lim_{\eta \to \lambda}e^\xi\eta$;
\item  \label{RecHypIV}for any $\xi$ and any limit $\lambda$, $e^\lambda(\xi+1) = \lim_{\eta \to \lambda}e^\eta(e^\lambda(\xi) + 1)$.
\end{enumerate}
\end{lemma}

\begin{cor}\label{cofexp}
If $\xi$ is a limit ordinal, then $\cf(e^\lambda\xi) = \cf(\xi)$. Otherwise,  $\cf(e^\lambda\xi) = \max(\cf(\lambda), \aleph_0)$.
\end{cor}



As a different application of Lemma \ref{recexp}, the reader might want to check that the function $\alpha \mapsto e^\omega(1 + \alpha)$ coincides with the epsilon function. Of course, Lemma \ref{recexp}.\ref{RecHypIV} is not the only way to compute $e^\lambda(\xi+1)$. Two useful alternatives are given as follows:





\begin{lemma}\label{explimits}
Let $\xi$ be an ordinal, $\lambda$ be an additively indecomposable limit, and $f$ be a nonzero-valued function on $\Ord$.
\begin{enumerate}[(i)]

\item\label{explimits1} If $0 < f\eta<e^\lambda (\xi+1)$ for all $\eta< \lambda$, then $\lim_{\eta \to \lambda}e^\eta(e^\lambda\xi + f\eta) = e^\lambda(\xi+1).$

\item\label{explimits2}
If $\gamma<\lambda$ is a limit and $0 < f\eta<e^\omega(e^\lambda\xi + 1)$ for all $\eta < \gamma$, then $\lim_{\eta \to \gamma}e^\eta(e^\lambda\xi + f\eta) = e^\gamma(e^\lambda\xi+1).
$
\end{enumerate}
\end{lemma}

\proof
For the first claim, we clearly have
\[\lim_{\eta \to \lambda}e^\eta(e^\lambda\xi + f\eta) 
\geq \lim_{\eta \to \lambda}e^\eta(e^\lambda\xi + 1)
= e^\lambda(\xi+1).\]
For the other inequality, fix $\eta<\lambda$ and take $\lambda_\eta \in(0, \lambda)$ such that $f\eta<e^{\lambda_\eta}(e^\lambda \xi+1)$. Since $e^{\lambda_\eta}(e^\lambda \xi+1)$ is additively indecomposable, we have that $e^\lambda\xi + f\eta < e^{\lambda_\eta}(e^\lambda \xi+1),$ whence $
e^\eta(e^\lambda\xi + f\eta) 
\leq e^{ \eta + \lambda_\eta  }(e^\lambda\xi + 1)
\leq \lim_{\zeta \to \lambda}e^{\zeta}(e^\lambda\xi + 1)
= e^\lambda(\xi+1).
$ Since $\eta$ was arbitrary, the claim follows.

For the second claim we have that, if $\eta < \gamma$, then $
e^\eta(e^\lambda\xi + f\eta) 
\geq e^{\eta}(e^{\lambda}\xi+1),$ and since $\lim_{\eta \to \gamma}e^{\eta}(e^{\lambda}\xi+1)
= \lim_{\eta \to \gamma}e^{\eta}(e^{\gamma}e^{\lambda}\xi+1)
= e^{\gamma}(e^{\lambda}\xi+1),$ it follows that $\lim_{\eta \to \gamma}e^\eta(e^\lambda\xi + f\eta) \geq e^\gamma(e^\lambda\xi+1)$. In addition, for any $\eta < \gamma$ and any $n < \omega$ such that $f\eta<e^{n}(e^{\lambda}\xi+1)$, we have that $e^\eta(e^\lambda\xi + f\eta) 
\leq e^{\eta+n}(e^{\lambda}\xi+1)
< e^{\gamma}(e^{\lambda}\xi+1),$ so that equality holds.
\endproof

\begin{lemma}[Hyperexponential normal form]\label{whnf}
For every ordinal $\xi>0$, there exist unique ordinals $\alpha$, $\beta$ such that $\beta$ is $1$ or additively decomposable and $\xi=e^\alpha \beta$.
\end{lemma}

\proof
Let $\alpha$ be the supremum of the set $A=\{\zeta: e^\zeta 1\leq \xi\}$. Note that this collection is indeed a set since it is a subset of $\xi+1$, while it is non-empty since $0\in A$. Moreover, $\alpha$ belongs to $A$, since the function $\zeta\mapsto e^\zeta 1$ is normal and thus $A$ is closed.

It follows that $\xi=e^\alpha \beta$ for some unique $\beta$. But $\beta$ is not in the range of $e$, since if we had $\beta=e\gamma$ then $\xi=e^{\alpha+1}\gamma$, contradicting the assumption that $\alpha$ is an upper bound for $A$. It follows from the definition of $e$ and Lemma \ref{indecomposable}.\ref{indecomposableOne} that $\beta$ is either $1$ or additively decomposable, as claimed.
\endproof

We call $\alpha$ above the {\em degree of indecomposability} of $\xi$; in particular, if $\xi$ is already additively decomposable, then $\alpha=0$. Note that by writing $\beta$ as a sum of indecomposables we may iterate this lemma and thus write any ordinal in terms of $e,+,0$ and $1$. This form is unique if we do not allow sums of the form $\xi+\eta$ for $\xi<\eta$. We next turn our attention to {\em hyperlogarithms,} which are iterations of {\em initial,} rather than normal, functions.

\begin{defi}[Initial function]
A function $f$ on the ordinals is said to be \emph{initial} if it maps initial segments to initial segments.
\end{defi}

\begin{defi}[Hyperlogarithms]
We define the \emph{hyperlogarithms} $(\ell^\xi)_{\xi \in \Ord}$ to be the unique family of initial functions that satisfy:
\begin{enumerate}[(i)]
\item\label{DefCoHI} $\ell^1 = \ell$,
\item\label{DefCoHII} $\ell^{\alpha + \beta} = \ell^\beta  \ell^\alpha$ and
\item if $(f^\xi)_{\xi\in\Ord}$ is a family of functions satisfying \ref{DefCoHI} and \ref{DefCoHII}, then for all $\alpha,\beta\in\Ord$, $\ell^\alpha\beta\geq f^\alpha\beta$.
\end{enumerate}
\end{defi}

As can be seen from the definition, hyperlogarithms are not right-additive but rather left-additive. This is to allow for surjectivity and the possibility that they have a right-inverse. Note, however, that they are not invertible since they are not injective (for example, $\ell 1=0=\ell (\omega + 1)$); we will, however, be interested in taking preimages under hyperlogarithms, in which case for a set of ordinals $A$ we write $\ell^{-\xi}(A)$ instead of $(\ell^{\xi})^{-1}(A)$.

Hyperlogarithms have some nice properties (see \cite{hyperations}):

\begin{lemma}[Properties of hyperlogarithms]\label{LemmPropLog}\
\begin{enumerate}[(i)]
\item $\ell\xi\leq\xi$ for all $\xi$.
\item\label{LemmPropLogII} If $\xi$ and $\delta$ are nonzero, then $\ell^\xi(\gamma + \delta) = \ell^\xi\delta$.
\item\label{LemmPropLogIII} For any $\gamma$, the sequence $(\ell^\xi\gamma)_{\xi\in \Ord}$ is non-increasing.
\item\label{LemmPropLogIV} If  $\lambda=\alpha+\omega ^\beta$ is a limit ordinal and $\xi$ is any ordinal, then there exists $\sigma<\lambda$ such that $\ell^\zeta \xi=e^{\omega^\beta}\ell^\lambda\xi$ for all $\zeta\in[\sigma,\lambda)$.
\end{enumerate}
\end{lemma}


Part \ref{LemmPropLogIII} of the above lemma states that hyperlogarithms are non-increasing in the exponent. They are not, however, monotone in their arguments: $\ell\omega = 1$ while $\ell(\omega+1)=0$. 


It is easy to see that $\ell$ functions as a left inverse for $e$. Moreover, \cite[Theorem 9.1]{hyperations} implies that $\ell^\xi$ is also a left inverse for $e^\xi$ for all $\xi$. Hence, if $\xi < \zeta$, then by rewriting $e^\zeta$ as $e^\xi e^{-\xi + \zeta}$, we obtain the following:

\begin{lemma} \label{exponentialcancellation}
If $\xi < \zeta$, then $\ell^\xi e^\zeta = e^{-\xi + \zeta}$ and $\ell^\zeta e^\xi = \ell^{-\xi + \zeta}$. Furthermore, if $\alpha < e^\xi\beta$, then $\ell^\xi \alpha < \beta$, and if $\alpha <\ell^\xi  \beta$, then $e^\xi\alpha <  \beta$.
\end{lemma}

We conclude this section with two examples, with aims at illustrating the properties and recursive computations of hyperexponentials and hyperlogarithms:

\begin{exm}
We calculate $e^{\omega_1}1$:
\begin{equation*}
e^{\omega_1}1 = 
\lim_{\eta \to {\omega_1}}e^\eta(e^{\omega_1}0 + 1) =
\lim_{\eta \to {\omega_1}}e^\eta(1).
\end{equation*}
On the one hand, each term $e^\eta1$ is smaller than $\omega_1$; on the other hand, for any $\eta<\omega_1$, we have $\eta\leq e^\eta 1 $. Therefore, the sequence must converge to $\omega_1$, whence $e^{\omega_1}1= \omega_1$.
\end{exm}

\begin{exm}
Consider the ordinal $\gamma = \omega_{\varepsilon_0} + \varepsilon_{\varepsilon_{\omega_1} + \varepsilon_{\omega^2\cdot 3}}$. We compute the sequence $(\ell^\xi\gamma)_{\xi \in \Ord}$.

Recall that $e^\omega$ equals the epsilon function. Thus, for $0<\xi<\omega$: 
\begin{equation*}
\ell^\xi\gamma = 
\ell^\xi e^\omega(\varepsilon_{\omega_1} + \varepsilon_{\omega^2 \cdot 3}) = 
e^{-\xi + \omega}(\varepsilon_{\omega_1} + \varepsilon_{\omega^2 \cdot 3}) = 
e^{\omega}(\varepsilon_{\omega_1} + \varepsilon_{\omega^2 \cdot 3}) = 
\varepsilon_{\varepsilon_{\omega_1} + \varepsilon_{\omega^2 \cdot 3}}.
\end{equation*}
When $\xi = \omega$, the exponential and the logarithm are cancelled: $\ell^\xi\gamma = \varepsilon_{\omega_1} + \varepsilon_{\omega^2 \cdot 3}$. For $\omega < \xi < \omega + \omega$, with $0<\zeta$, we have 
\begin{equation*}
\ell^\xi\gamma = 
\ell^{-\omega + \xi} (\varepsilon_{\omega_1} + \varepsilon_{\omega^2 \cdot 3}) =
\ell^{-\omega + \xi} e^\omega (\omega^2 \cdot 3) = 
\varepsilon_{\omega^2 \cdot 3}.
\end{equation*}
As above, if $\xi = \omega \cdot 2$, then $\ell^\xi\gamma = \omega^2 \cdot 3$. $\omega^2\cdot 3$ can be written as $\omega^2\cdot 2 + \omega^2$, which implies that $\ell^{\omega \cdot 2 + 1} = \ell\omega^2 = \ell (e2) = 2$ and $\ell^{\omega\cdot 2 + 2} = \ell2 = 0$. For any $\xi$ greater than $\omega\cdot 2 + 2$, we also have $\ell^\xi\gamma = 0$.
\end{exm}

Hyperexponential ordinal notations will be crucial in the description of topological semantics in later sections, while hyperlogarithms will play an important role in the description of the topologies that we will use. 

\section{Topological interpretation of provability logic}\label{SecTopo}

Recall that a \emph{topological space} is a pair $(X, \tau)$, where $X$ is a set and $\tau$ is a family of subsets of $X$ containing $X$ and $\varnothing$ that is closed under arbitrary unions and finite intersections. An open set $U$ containing a point $x$ is a \emph{neighborhood} of $x$. The set $U \setminus \{x\}$ is then called a \emph{punctured neighborhood} of $x$. For any $A \subset X$, we say that $x$ is a \emph{limit point} of $A$ if $A$ intersects every punctured neighborhood of $x$. We call the set of limit points of $A$ the \emph{derived set} of $A$ and denote it by $dA$. We will also denote it by $d_\tau A$ if we want to emphasize the topology with respect to which we take the limit points. A point in $A$ that is not a limit point is \emph{isolated} in $A$. Equivalently, a point $x$ is isolated in $A$ if, and only if, $\{x\}$ is open in $A$ under the inherited topology (recall that for any topological space $X$ and any subset $S$, the \emph{inherited topology} is formed by all sets $U \cap S$, where $U \in \tau$). We can define topological semantics for the language $\mlang$ by assigning subsets of $X$ to each propositional variable, in the following way:

\begin{defi}[Topological semantics]\label{DefKripkSem}
Let $\mathfrak{X}=(X, \tau)$ be a topological space. A \emph{valuation} is a function $\lb\cdot\rb:L_\Box \to \wp(X)$ such that for any $\varphi, \psi$, $\lb\bot\rb = \varnothing$, $\lb\lnot \varphi\rb = X \backslash \lb\varphi\rb$, $\lb\varphi \land \psi\rb = \lb\varphi\rb \cap \lb\psi\rb$ and $\lb\Diamond \varphi\rb = d \lb\varphi\rb$.

A \emph{model} $\mathfrak{M} = (\mathfrak{X}, \lb\cdot\rb)$ is a topological space together with a valuation. If $\mathfrak{M}$ is a model, we may write $\mathfrak{M},x \models \varphi$ instead of $x \in \lb\varphi\rb$; in this case, we say $\varphi$ is {\em true} at $x$ (in $\mathfrak{M}$). We say that $\varphi$ is \emph{satisfied} in $\mathfrak{M}$ if $\lb\varphi\rb$ is nonempty, that $\varphi$ is \emph{true} in $\mathfrak{M}$ and write $\mathfrak{M} \models \varphi$ if $\lb\varphi\rb = W$, and we say $\varphi$ is \emph{valid} in a space $\mathfrak{X}$ and write $\mathfrak{X} \models \varphi$ if $\lb\varphi\rb = X$ for any model $(\mathfrak{X},\lb\cdot\rb)$. 
\end{defi}

Interpreting the modal diamond as the derived set operator is sometimes called the \emph{d-interpretation} of modal logic \cite{bezhanishvilimorandi}. An important example is given as follows:

\begin{defi}[Relational structure]\label{example trees}
For our purposes, a \emph{relational structure} is a set $T$ together with a transitive, irreflexive relation $R$ on $T$. Given a relational structure $(T, R)$, we assign a topology $\tau_R$, called the {\em upset topology,} to $T$, whose open sets are those $U\subset W$ such that, whenever $x\in U$ and $xRy$, it follows that $y\in U$.
\end{defi}

For those familiar with Kripke semantics, it is straightforward to check that for these structures, the usual Kripke interpretation based on $R$ and the $d$-interpretation based on $\tau_R$ coincide. 

\begin{defi}[Soundness and completeness]
Let $\mathcal{X}$ be a class of topological spaces and $\sf L$ be a logic over $\mlang$. We say $\sf L$ is \emph{sound} with respect to $\mathcal{X}$ if $\mathfrak X \models \varphi$ for all $\mathfrak X \in \mathcal{X}$ whenever $\mathsf L \vdash \varphi$. Conversely, we say $\mathsf L$ is \emph{complete} with respect to $\mathcal{X}$ if whenever $\mathfrak X \models \varphi$ for all $\mathfrak X \in \mathcal{X}$, then $\mathsf L \vdash \varphi$. Moreover, $\mathsf L$ is \emph{strongly complete} with respect to $\mathcal{X}$ if for any set $\Gamma$ of $\mlang$-formulae that is consistent with respect to $\mathsf L$, there is some $\mathfrak X \in \mathcal{X}$ such that $\Gamma$ is satisfied in $\mathfrak X$.
\end{defi}

It is easy to verify that strong completeness implies simple completeness.

\begin{defi}[Scattered space]
We say that a topological space $(X,\tau)$ is \emph{scattered} if each nonempty subset of $X$ has an isolated point.
\end{defi}

\begin{exm}
Recall that a relation $R$ on a set $T$ is \emph{converse well-founded} if every nonempty subset of $T$ has an $R$-maximal element. Then, if $R$ is converse well-founded, $(T, \tau_R)$ is a scattered space.
\end{exm}





\begin{lemma}[\cite{esakia1981}]
L\"ob's axiom is valid in a topological space $\mathfrak X$ if and only if $\mathfrak X$ is scattered.
\end{lemma}

It is straightforward to check that the axiom $K$ is valid and the rules of necessitation and modus ponens preserve validity in all topological spaces (see, for example, \cite{jvblogicsspace}). With this and the preceding discussion we may formulate Esakia's theorem \cite{esakia1981}:

\begin{theorem}
A formula $\varphi$ is a theorem of $\gl$ if, and only if, it is valid in all scattered spaces.
\end{theorem}

There are several possible improvements to Esakia's theorem. Recall that a \emph{tree} is a relational structure $(T,R)$ (which by our definition is transitive), such that each $R$-predecessor set is well-ordered, and $T$ contains a unique $R$-minimal element, its {\em root.} The following was originally proved by K. Segerberg \cite{segerberg1971}:

\begin{theorem}\label{treecompletenessglp1}
$\gl$ is sound and complete with respect to the class of finite trees.
\end{theorem}



Nonetheless, $\gl$ is not strongly complete with respect to this semantics, as is well known:

\begin{exm}
Consider the set $\Gamma = \{\Diamond p_0\} \cup \{\Box(p_i \rightarrow \Diamond p_{i+1}) \colon i < \omega\}.$ It is easy to see that any finite subset of $\Gamma$ is satisfiable in a finite tree, whence $\Gamma$ is consistent with $\gl$ by soundness. Moreover, all of $\Gamma$ is not satisfiable in any converse well-founded tree. We leave the details to the reader.
\end{exm}

Some more examples of scattered spaces are provided by natural topologies on ordinal numbers.

\begin{exm}\ \label{examples of scattered spaces}
\begin{enumerate}
\item Let $\Theta$ be an ordinal number and $\tau$ consist of all sets $[0,\beta]$. Then $(\Theta, \tau)$ is a scattered space: the least element of each set is isolated. We call this topology the \emph{left topology} and denote it $\mathcal{I}_0$. Whenever an ordinal is equipped with the left topology, we may write $\Theta_0$ instead of $(\Theta, \mathcal{I}_0)$.

\item Let $\Theta$ be an ordinal number and $\tau$ consist of all sets $(\alpha,\beta]$ and all sets $[0,\beta]$. Then $(\Theta, \tau)$ is a scattered space: each successor ordinal is isolated. We call this topology the \emph{interval topology} and denote it $\mathcal{I}_1$. Whenever an ordinal is equipped with the interval topology, we may write $\Theta_1$ instead of $(\Theta, \mathcal{I}_1)$.

\item Let $\Theta$ be an ordinal number. We define the \emph{club topology}, $\tau_c$ to be the unique topology such that $U$ contains a neighborhood of $\xi \in U$ if, and only if, $U$ contains a club (i.e., $\ico 1$-closed and unbounded set) in $\xi$ or $\cf(\xi) < \aleph_1$. $(\Theta, \tau_c)$ is a scattered space, for if $\aleph_1 \leq \cf(\xi)$ and $A$ is a club on $\xi$, then $A$ contains a point $\zeta$ such that $\cf(\zeta) < \aleph_1$.
\end{enumerate}
\end{exm}

In the following sections, we will introduce further improvements to Esakia's theorem and in fact exhibit single topological spaces with respect to which $\gl$ is {\em strongly} complete. In order to do this, we must study scattered spaces in more detail. If $(X,\tau)$ is a topological space, we define transfinite iterations $( d^\xi )_{\xi\in\Ord}$ of the derived-set operator by $d^0A = A$, $d^{\xi+1}A = dd^\xi A$ and $d^\lambda A = \bigcap_{\xi<\lambda}d^\xi A$ for limit $\lambda$.

The derived set operator and its iterations have the following properties, which are easy to prove (see, for example, \cite{bezhanishvilimorandi}):

\begin{lemma} \label{properties of derivatives}Let $(X,\tau)$ be a scattered space.
\begin{enumerate}[(i)]
\item If $A$ is closed in $X$, then $dA$ is closed in $X$;
\item \label{DerClosed}$d^\alpha X$ is closed for each $\alpha$;
\item \label{DerInclusion}if $\alpha \leq \beta$, then $d^\beta X \subset d^\alpha X$.
\end{enumerate}
\end{lemma}

By Lemma \ref{properties of derivatives}.\ref{DerInclusion}, iterations of the derived-set operator stabilize at some stage below the successor of $\vert X \vert$. The following was famously noted by Cantor:
\begin{prop}
A topological space $(X,\tau)$ is scattered if and only if $d^\xi X =\varnothing$ for some $\xi$.
\end{prop}

This leads naturally to a notion of the ``size'' of a scattered space--- its rank:

\begin{defi}[Rank]
Let $\mathfrak X=(X,\tau)$ be a scattered space. For any $x \in X$, we define $\rho_\tau x$, the \emph{rank} of $x$, to be the least ordinal $\xi$ such that $x \not \in d^{\xi+1} X$. We may also write it as $\rho_\mathfrak X x$, or even $\rho_X x$ or $\rho x$ if there is no risk of confusion. We also define the rank of $\mathfrak X$ to be $\rho \mathfrak X = \sup_{x \in X}(\rho_\tau x + 1)$. We may also write it as $\rho_\tau X$ or even $\rho X$ when there is no risk of confusion.
\end{defi}

One useful fact about the rank function is that near any point $x$, we can always find points of every smaller rank:

\begin{lemma} \label{propertiesofdmapsTwo} 
If $V$ is a neighborhood of $x$ in a scattered space, then $[0, \rho x) \subset \rho(V \setminus \{x\})$. Furthermore, there is a neighborhood of $x$ for which equality holds.
\end{lemma}

As way of illustration, let us list the rank functions in some important scattered spaces.

\begin{prop}\
\begin{enumerate}
\item \label{rankInit} The rank function on any ordinal space $\Theta_0$ is the identity: $\rho_0(\xi) = \xi$. 
\item \label{rankInt} The rank function on $\Theta_1$ is the end-logarithm: $\rho_1(\xi) = \ell\xi$. 
\item \label{rankclub} The rank function on $(\Theta, \tau_c)$ is given as follows: let $(\Regular\alpha)_{\alpha\in\Ord}$ enumerate all infinite regular cardinals. Then, for nonzero $\alpha$, $\rho_{\tau_c}(\xi)=\alpha$ if, and only if, $\cf(\xi)=\Omega_\alpha$.
\end{enumerate}
\end{prop}

\proof
Items \ref{rankInit} and \ref{rankInt} are instances of Lemma \ref{LemmHyperRank} below. To prove item \ref{rankclub}, we show by induction that if $0 < \alpha$, then $\xi\in d^\alpha_{\tau_c}\Theta$ if and only if $\cf(\xi)\geq \Regular\alpha$. We proceed by induction on $\alpha$. The claim is vacuously true when $\alpha=0$, and readily follows from the induction hypothesis when $\alpha$ is a limit since in this case $d^\alpha_{\tau_c}\Theta=\bigcap_{\beta<\alpha}d^\beta_{\tau_c}\Theta$.

Thus we may assume $\alpha = \beta + 1$ and $\xi \in d^\beta_{\tau_c}\Theta$. By induction hypothesis, $\Regular\beta \leq \cf(\xi)$. If the inequality is strict, then $\xi$ is a limit point of $d^\beta_{\tau_c}\Theta$, as is it easy to check that $\{\zeta < \xi \colon \cf(\zeta) = \Regular\beta\}$ is stationary in $\xi$\jp{Viene en Jech (p. 94, aunque s\'olo para cardinales regulares, el argumento es el mismo), pero no viene con la prueba. S\'olo dice "es f\'acil verlo".} (i.e., it intersects every club). By \cite[Theorem 11.5]{bekgab14}, $\tau_c$ is generated by $\ico 1$ and sets $\{d_{\ico 1}A \colon A \subset \Theta\}$, so that if, on the contrary, $\Regular\beta = \cf(\xi)$, then $\{\xi\} = d^\beta_{\tau_c} \Theta \cap d_{\ico 1}A$, where $A$ is any set cofinal in $\xi$ of order-type $\Regular\beta$ and so $\xi$ is isolated in $d^\beta_{\tau_c}\Theta$.
\endproof

It is crucial to define appropriate structure-preserving mappings between scattered spaces. Obviously, homeomorphisms preserve all the relevant structure but, as it turns out, a weaker condition will do for our purposes. Below, a function between topological spaces is \emph{pointwise discrete} if the preimage of any singleton is a discrete subspace.

\begin{defi}[$d$-map]
Let $X$ and $Y$ be scattered spaces. Then $f:X \to Y$ is a \emph{$d$-map} if it is continuous, open, and pointwise discrete.
\end{defi}

An important example is one that we introduced above, as shown in \cite{topocompletenessofglp}:

\begin{lemma}\label{uniquedmap}
For any scattered space $\mathfrak X = (X,\tau)$ and any ordinal $\Theta\geq \rho \mathfrak X$, the rank function is the unique $d$-map $f:\mathfrak X \to \Theta_0$.
\end{lemma}




Obviously, any homeomorphism is a $d$-map and so, since the composition of $d$-maps is a $d$-map, they can be thought of as morphisms in the category of scattered spaces. Two key facts about $d$-maps is that they preserve ranks and the validity of formulae between topological spaces if they are surjective (see \cite{topocompletenessofglp} for item \ref{propertiesofdmapsOne} and \cite{bmm} for item \ref{valuation transfer}).

\begin{lemma}\ \label{propertiesofdmaps}
\begin{enumerate}[(i)]
\item\label{propertiesofdmapsOne} If $f:\mathfrak X \to \mathfrak Y$ is a $d$-map, then $\rho_\mathfrak X = \rho_\mathfrak Y f$.
\item\label{valuation transfer} Suppose $\mathfrak X$ and $\mathfrak Y$ are scattered spaces such that a surjective $d$-map $f: \mathfrak X \to \mathfrak Y$ exists. Then for any $\mlang$-formula $\varphi$, $\mathfrak X \models \varphi$ implies $\mathfrak Y \models \varphi$.
\end{enumerate}
\end{lemma}

\section{$\omega$-bouquets}\label{SecRel}

As we have seen, $\gl$ is not strongly Kripke-complete. However, it turns out that this situation may be remedied with only a minor modification of its relational semantics.

\begin{defi}[$\omega$-bouquet]\label{omegabouquets}
Let $(T,R)$ be a countable, converse well-founded tree, and let $\rho\colon T\to \Ord$ be the rank function on $T$ with respect to the upset topology.

We define a new topology, $\sigma_R$, to be the least topology extending $\tau_R$ such that if $w\in T$ is such that $\rho ( w )$ is a limit ordinal, $\{v_i\}_{i<\omega}$ enumerates all daughters of $w$ exactly once, and $n<\omega$, then $\{w\}\cup \bigcup _{n<i} \big ( \{v_i\}\cup R(v_i) \big ) \in \sigma_R$.

We say a topological space $(T,\sigma)$ is an {\em $\omega$-bouquet} if there exists a binary relation $R$ on $T$ such that $(T,R)$ is a countable, converse well-founded tree and $\sigma=\sigma_R$.
\end{defi}

Thus, points of limit rank have new punctured neighbourhoods, consisting of the union of upsets of all but finitely many immediate successors. An easy exercise is to check that $\sigma_R$ is a rank-preserving extension of $\tau_R$. This apparently innocuous change will give us a class of structures for which $\gl$ is strongly complete. A similar construction was used by V. Shehtman \cite{shehtman99} to prove the strong completeness of many modal logics for their topological semantics, from which we adopt the term {\em bouquets}.

In this section, we will prove that $\gl$ is strongly complete for the class of $\omega$-bouquets (i.e., any consistent set of formulae $\Gamma$ is satisfiable on an $\omega$-bouquet). Without loss of generality, we will assume $\Gamma$ to be maximal consistent. Recall that we are asuming that the set of propositional variables is countable, so that $\Gamma$ will be countable as well. We introduce some auxiliary notation:

\begin{defi}
If $\Gamma$ is a set of $\gl$ formulae, define $\Gamma^\Box=\{\varphi:\Box\varphi\in \Gamma\}$ and $\Box\Gamma=\{\Box\varphi:\varphi\in\Gamma\}$.
\end{defi}

\begin{defi}
If $\Gamma$ is a maximally consistent set of $\gl$ formulae, we define the {\em characteristic} of $\Gamma$ as the supremum of $\{n < \omega : \Diamond^n \top \in \Gamma\}$.
\end{defi}

Note that the characteristic of $\Gamma$ may be either finite or $\omega$. In our proof, we will consider each of these two cases separately. First, we will consider the case where $\Gamma$ has finite characteristic. We will need the following lemma:

\begin{lemma}
Suppose that $\Gamma$ is maximally consistent with characteristic $n+1<\omega$, and $\Diamond \varphi\in \Gamma$. Then, $\Delta=\{\varphi,\Box^{n+1}\bot\}\cup \Gamma^\Box\cup \Box\Gamma^\Box$ is consistent.
\end{lemma}

\proof
Towards a contradiction, assume that $\Delta$ is inconsistent, so that there are $\psi_0,...,\psi_{m-1}\in \Gamma^\Box$ such that
\[\vdash \Box^{n+1}\bot\wedge\bigwedge_{i<m} (\psi_i\wedge\Box\psi_i) \to \neg\varphi.\]

Reasoning in $\sf GL$, it follows that
\[\vdash \Box^{n+2}\bot\wedge\bigwedge_{i<m} \Box(\psi_i\wedge\Box\psi_i) \to \Box\neg\varphi.\]
But $\Box(\psi_i\wedge\Box\psi_i)$ is equivalent to $\Box\psi_i$, so that
\[\vdash \Box^{n+2}\bot\wedge\bigwedge_{i<m} \Box\psi_i \to \Box\neg\varphi.\]
Since $\Gamma$ has characteristic $n+1$ and is maximal consistent, we have that $\Box^{n+2}\bot\in\Gamma$, while also each $\Box\psi_i\in\Gamma$. Thus $\Gamma\vdash \Box\neg\varphi$, together with $\Diamond\varphi\in\Gamma$, contradicts the consistency of $\Gamma$.
\endproof

As an immediate consequence, we obtain:

\begin{lemma}\label{LemmPhiGamm}
Suppose that $\Gamma$ is maximally consistent with characteristic $n+1<\omega$, and $\Diamond \varphi\in \Gamma$. Then, there is a maximal consistent set $\Gamma_\varphi$ of some characteristic $m\leq n$ such that $\Gamma_\varphi\supset\{\varphi\}\cup \Gamma^\Box\cup \Box\Gamma^\Box$.
\end{lemma}

With this, we may prove our strong consistency result in the case of finite characteristic.

\begin{lemma}\label{LemmFinCh}
Let $\Gamma$ be maximal consistent and of finite characteristic $n$. Then, $\Gamma$ is satisfiable in a model based on an $\omega$-bouquet of rank $n+1$.
\end{lemma}

\proof
By induction on $n$. Assume by induction hypothesis that the lemma holds for any $\Gamma'$ with characteristic $m<n$. Thus, for each $\Diamond\varphi\in \Gamma$, there is a model $(T_\varphi,R_\varphi,\lb\cdot\rb_\varphi)$ of rank $n_\varphi<n$ satisfying $\Gamma_\varphi$, defined as in Lemma \ref{LemmPhiGamm}. Moreover, $n_\varphi=n-1$ if $\varphi=\Diamond^{n-1}\top\in\Gamma$.

Define $(T,R)$ by taking the disjoint union of all $(T_\varphi,R_\varphi)$ and adding a fresh root $r$. We consider a valuation $\lb \cdot \rb$ over $T$ induced by the valuations $\lb\cdot\rb_\varphi$ described as follows. Define $\lb \cdot \rb_0$ by setting $r \in \lb p \rb_0$ if and only if $p\in\Gamma$, and then set $\lb p \rb=\lb p \rb_0\cup\bigcup_{\Diamond\varphi\in\Gamma}\lb p \rb_\varphi$ for each propositional variable $p$. It is then straightforward to check that $r\in\lb \varphi\rb$ if and only if $\varphi\in\Gamma$.
\endproof

It remains to exhibit models that satisfy sets of characteristic $\omega$. We will do so by ``decomposing'' $\Gamma$ into smaller pieces and applying Theorem \ref{treecompletenessglp1} to each of them.

\begin{defi}
Assume that $\Gamma$ has characteristic $\omega$. Let $(\psi_i)_{i < \omega}$ enumerate all formulae such that $\Diamond\psi_i\in\Gamma$, in a way such that each $\psi_i$ occurs infinitely often; let
$(\phi_i)_{i < \omega}$ enumerate all formulae such that $\Box\phi_i\in\Gamma$ (in any way). For $i<\omega$, define $\Gamma(i)=\{ \psi_{i} \} \cup \{ \phi_j\wedge \Box\phi_j\}_{j<i}.$
\end{defi}

\begin{lemma}\label{LemmInfChar}
If $\Gamma$ is consistent, then so is each $\Gamma(i)$. Moreover, each $\Gamma(i)$ is satisfied in a finite model $(T_i,R_i,\lb\cdot\rb_i)$.
\end{lemma}

\proof
It is easy to see that $\Gamma\vdash\Diamond\bigwedge\Gamma(i)$ for all $i$, so the latter must be consistent. Since $\Gamma(i)$ is finite, we can apply Theorem \ref{treecompletenessglp1}.
\endproof

With this, we are ready to prove the main result of this section.

\begin{theorem}\label{TheoCompBouq}
$\gl$ is strongly complete for the class of $\omega$-bouquets.
\end{theorem}

\proof
Let $\Gamma$ be a consistent set of formulae; without loss of generality, assume $\Gamma$ to be maximal. If $\Gamma$ has finite rank, then apply Lemma \ref{LemmFinCh}.

Otherwise, $\Gamma$ has rank $\omega$. By Lemma \ref{LemmInfChar}, each $\Gamma(i)$ is satisfied in a finite model $(T_i,R_i,\lb\cdot\rb_i)$.

As in the proof of Lemma \ref{LemmFinCh}, take the disjoint union of all $(T_i,R_i,\lb\cdot\rb_i)$, and add a root $r$ so that, for all variables $p$, $r \in \lb p\rb$ if and only if $p\in\Gamma$; call the resulting model $\mathfrak M$. It is then easy to check that $\mathfrak M,r\models \Gamma$.
\endproof

Theorem \ref{TheoCompBouq} is interesting on its own right, but our focus is on ordinal spaces. Our strategy will be to `lift' this result using $d$-maps through Lemma \ref{propertiesofdmaps}.\ref{valuation transfer}. Observe that our construction uses only very specific $\omega$-bouquets; they are either of finite rank, or of rank $\omega+1$, consisting of a root added to countably many finite subtrees. However, it will be just as easy to construct $d$-maps onto an arbitrary $\omega$-bouquet as onto one of the above forms.

Let us move on to define the main scattered spaces in which we will be interested.

\section{Generalized Icard spaces}\label{SecIc}
In this section, we will discuss \emph{generalized Icard spaces based on a scattered space $\mathfrak X$,} which are an increasing sequence of topologies $(\ic\lambda{\mathfrak X})_{\lambda\in\Ord}$. These topologies were introduced with $\mathfrak X$ of the form $(\Theta,\mathcal I_0)$ to give semantics for the variable-free fragment of $\glp_\omega$ in \cite{icardglp}. They were then defined for arbitrary $\glp_\Lambda$ in \cite{modelsofglp} and subsequently studied in \cite{polytopologies}. We generalize them further by letting $\mathfrak X$ be an arbitrary scattered space. As a notational convention, we will sometimes write ordinal intervals $[0, \alpha]$ as $(-1, \alpha]$.

\begin{defi}[Icard topology]
Let $\mathfrak X=(X, \tau)$ be a scattered space of rank $\Theta$. Given an ordinal $\lambda$, we define a topology $\ic \lambda\tau$ to be the least topology containing $\tau$ and all sets of the form
\begin{equation*}
(\alpha, \beta]_\xi^{\mathfrak X} = \{x \in X: \alpha < \ell^\xi \rho_\tau x \leq \beta\},
\end{equation*}
for some $-1 \leq \alpha < \beta \leq \Theta$ and some $\xi < \lambda$. These are called the \emph{generalized Icard topologies based on $\mathfrak X$}. We call spaces of the form $\ic\lambda{\mathfrak X}=(X,\ic \lambda\tau)$ for some $\tau$ \emph{generalized Icard spaces}. 
\end{defi}

If $\mathfrak X=(X,\tau)$, we may write $(\alpha, \beta]_\xi^{\tau}$ instead of $(\alpha, \beta]_\xi^{\mathfrak X}$. We may also index properties by an ordinal $\lambda$ if they refer to the topology $\ic \lambda\tau$. In this way, $\rho_{\ic\lambda\tau}$ becomes the $\lambda$-rank $\rho_\lambda$.

In the case of ordinal spaces, hyperlogarithms allow us to compute $\lambda$-ranks. Recall that we denote the left topology on the ordinals by $\mathcal {I}_0$. Thus, $\ic 0{(\mathcal{I}_0)}=\mathcal{I}_{0}$ and $\ic 1{(\mathcal{I}_0)}=\mathcal{I}_{1}$. In this case, we will write simply $\ico \lambda$ instead of $\ic \lambda{\ico 0}$ and define $\Theta_\lambda=(\Theta,\mathcal I_\lambda)$. The following is proved in \cite{polytopologies}:

\begin{lemma}\label{LemmHyperRank}
If $\Theta$, $\lambda$, and $\xi < \Theta$ are ordinals, 
then $\rho_\lambda\xi=\ell^\lambda\xi$. Moreover, if $\mu$ is any ordinal, then $\ell^\lambda: \Theta_{\lambda+\mu} \to \Theta_\mu$ is a $d$-map
\end{lemma}

Let us now prove some additional properties of generalized Icard spaces:

\begin{lemma}\label{nextlambdamap}
Let $\mathfrak X$ and $\mathfrak Y$ be scattered spaces and $f:\mathfrak X\to \mathfrak Y$ be a $d$-map. Then, $f:\ic\lambda{\mathfrak X} \to \ic\lambda{\mathfrak Y}$ is a $d$-map for any ordinal $\lambda$.
\end{lemma}

\proof
By Lemma \ref{propertiesofdmaps}.\ref{propertiesofdmapsOne}, $\rho_\mathfrak X=\rho_\mathfrak Y f$, so that for any $-1 \leq \alpha < \beta$, any $\xi < \lambda$, and all $\tau$-open $A$ and $\sigma$-open $B$, $f\big(A \cap (\alpha, \beta]_\xi^{\mathfrak X}) = f(A) \cap (\alpha, \beta]_\xi^{\mathfrak Y}$ and $f^{-1}(B) \cap (\alpha, \beta]_\xi^{\mathfrak X} = f^{-1}(B \cap (\alpha, \beta]_\xi^{\mathfrak Y}$). From this it follows that $f$ is open and continuous. Furthermore, $f$ is pointwise discrete because $\tau \subset \ic \lambda\tau$. 
\endproof

\begin{cor}\label{lambdadmaps}
Let $\mathfrak X$ be a scattered space of rank $\Theta$. Then for any ordinal $\lambda$, we have that $\rho_\mathfrak X:\ic\lambda{\mathfrak X} \to \Theta_\lambda$ is a surjective $d$-map.
\end{cor}

\proof Immediate from Lemmas \ref{uniquedmap} and \ref{nextlambdamap}.
\endproof

Corollary \ref{lambdadmaps} provides us with a nice way of determining the rank function of Icard topologies:
\begin{lemma}\label{icardlambdarank}
If $\mathfrak X$ is a scattered space, then  $\rho_{\lambda}=\ell^{\lambda}\rho_\mathfrak X$.
\end{lemma}

\proof
By Corollary \ref{lambdadmaps}, $\rho_\mathfrak X:\ic\lambda{\mathfrak X} \to (\rho \mathfrak X)_\lambda$ is a $d$-map, and the topology in $(\rho\mathfrak X)_\lambda$ is precisely $\mathcal{I}_\lambda$. Since, by lemma \ref{LemmHyperRank}, $\ell^\lambda$ is the rank function on $(\rho\mathfrak X)_\lambda$, we have that $\ell^{\lambda}\rho_\mathfrak X:\ic\lambda{\mathfrak X} \to (\rho \mathfrak X)_0$ is a $d$-map and coincides with the rank function by Lemma \ref{uniquedmap}. 
\endproof

Our choice of notation is explained by the following result:

\begin{lemma}\label{sums}
For any scattered space $(X, \tau)$, and any ordinals $\lambda$ and $\mu$:
\begin{equation*}
{\ic\mu{(\ic\lambda\tau)}}={\ic{(\lambda + \mu)}\tau}.
\end{equation*}
\end{lemma}

\proof
${\ic\mu{(\ic\lambda\tau)}}$ is the topology generated by $\ic\lambda\tau$ and all sets $(\alpha, \beta]^{\ic\lambda\tau}_\xi$, for $\xi < \mu$. $\ic\lambda\tau$ is generated by $\tau$ and all sets $(\alpha, \beta]^{\tau}_\xi$, for $\xi < \lambda$; and for any $\alpha$, $\beta$, and $\xi<\mu$:
\begin{align*}
(\alpha, \beta]^{\ic\lambda\tau}_\xi
&= \{x \in X: \alpha < \ell^{\xi}\rho_{\ic\lambda\tau}x \leq \beta\}= \{x \in X: \alpha < \ell^{\xi}\ell^{\lambda}\rho_\tau x \leq \beta\}\\
&= \{x \in X: \alpha < \ell^{\lambda + \xi}\rho_\tau x \leq \beta\}= (\alpha, \beta]^{\tau}_{\lambda + \xi},
\end{align*}
so that ${\ic\mu{(\ic\lambda\tau)}}$ is generated by $\tau$ and all sets $(\alpha, \beta]^{\tau}_{\xi}$, for $\xi < \lambda + \mu$, and hence is equal to $\ic{(\lambda+\mu)}\tau.$
\endproof

Another example of $d$-maps, which we present without proof, is given by ordinal subtraction:

\begin{lemma}\label{disphomeomorphism}
Let $\zeta,\lambda, \Theta $ be ordinals. Then, $(-\zeta + \cdot):[\zeta,\zeta+\Theta]_{\lambda} \to [0,\Theta]_{\lambda}$, i.e., the function given by $\xi \to -\zeta + \xi$, is a homeomorphism.
\end{lemma}

The following result shows that each point can be separated from points of equal rank in Icard topologies:

\begin{lemma}
Let $(X,\tau)$ be a scattered space and $\lambda$ be an ordinal. Any $x$ in $X$ has a $\lambda$-neighborhood $U$ such
that whenever $x \neq y \in U$, $\ell^\lambda \rho_0 y < \ell^\lambda \rho_0 x$. 
\end{lemma}

\proof
This is immediate from Lemma \ref{propertiesofdmapsTwo}, using the fact that $\rho_\lambda=\ell^\lambda\rho_0$ by Lemma \ref{icardlambdarank}.
\endproof

Recall that for any topological space $(X,\tau)$, a \emph{neighborhood base} for $x$ is a collection of open sets $\mathcal{N}$ such that for any open set $U$ containing $x$, there exists an open set $V \in \mathcal{N}$, also containing $x$, such that $V \subset U$. We will now provide neighborhood bases for Icard topologies in two useful ways. 

\begin{defi}
A \emph{simple function} is a function $r\colon A\to \Ord$ with $A \subset \Ord$ a finite set.

If $\mathfrak X$ is a scattered space, $r$ is a simple function and the rank of $x \in X$ is an ordinal $\alpha$ such that $r(\xi)<\ell^\xi\alpha$ for all $\xi \in \dom(r)$, define $B^\mathfrak X_r(x)$ to be the set  $\bigcap_{\xi \in dom(r)}(r(\xi),\ell^\xi\alpha]^\mathfrak X_\xi.$
\end{defi}

As usual, if $\mathfrak X=(X,\tau)$, we may write $B^\tau_r(x)$ instead of $B^\mathfrak X_r(x)$. Note that any set $B^\tau_r(\alpha)$ is an open set of $\ic\lambda\tau$ if $\max (\dom(r)) < \lambda$. In fact:

\begin{lemma}\label{nhbase1}
Let $\mathfrak X=(X,\tau)$ be a scattered space, $0 < \lambda$, $0 < \ell^\lambda\xi$, and $x \in X$ have rank $\xi$. Then the sets of the form $U \cap B_r^\tau(\xi)$ with $U\in \tau$ and $\dom(r) \subset \lambda$ form a neighborhood base for $x$ in $\ic\lambda\tau$. Moreover, if $\mathfrak X$ is an ordinal with the left topology, we can take $U = X$.
\end{lemma}

\proof
By definition, every neighborhood of $\xi$ in $\ic\lambda\tau$ contains another neighborhood of the form $V = U \cap \bigcap_{k < K} (\alpha_k, \beta_k]^\mathfrak X_{\sigma_k}$ for some $U \in \tau$, where all $\sigma_k < \lambda$.

We may assume that all $\alpha_k$ are greater than $-1$ since $0 < \ell^\lambda\xi \leq \ell^{\sigma_k}\xi$ for any $k$ and that all $\sigma_k$ are distinct because $\sigma_j = \sigma_k$ implies $(\alpha_k, \beta_k]^\mathfrak X_{\sigma_k} \cap (\alpha_j, \beta_j]^\mathfrak X_{\sigma_k} = (\max\{\alpha_k, \alpha_j\}, \min\{\beta_k,\beta_j\}]^\mathfrak X_{\sigma_k}$. To obtain the desired result, define $r$ to have as domain the set of all $\sigma_k$, so that $r(\sigma_k)=\alpha_k$ and note that $U\cap B_r^\mathfrak X(x) \subset V$.

Finally, if $\mathfrak X$ is an ordinal with the left topology, then by definition $(0, \xi]_0\subset U$. If we define $s$ so that $\dom(s)=\{0\}\cup \dom(r)$, $s(\eta)=r(\eta)$ when defined, and $s(0)=0$ if $r(0)$ is undefined, then we also obtain $X\cap B_s^\mathfrak X(x)=B_s^\mathfrak X(x) \subset V$, as claimed.
\endproof

\begin{lemma}\label{nhbase2}
Let $\mathfrak X=(X,\tau)$ be a scattered space, $\lambda=\alpha+\omega^\beta$ be an ordinal, and $x \in X$ have rank $e^\lambda \Theta$ with $\Theta>0$. Then, sets of the form $U \cap (\eta,e^{\omega^\beta}\Theta]^{\mathfrak X}_{\gamma}$, where $U \in \tau$, $\eta<e^{\omega^\beta} \Theta$, and $\gamma < \lambda$, form a $\ic \lambda\tau$-neighborhood base for $x$. If $\mathfrak X$ is an ordinal with the left topology, then we may take $U=[0,e^\lambda\Theta]_0$.
\end{lemma}

\proof Let $V$ be a $\ic\lambda\tau$-neighborhood of $x$. By Lemma \ref{nhbase1}, $V$ contains a set of the form $U \cap B_r^\tau(e^\lambda \Theta)$ for some simple function $r$ with domain contained in $\lambda$. We may assume using Lemma \ref{propertiesofdmapsTwo} that $U\setminus \{x\}\subset [0,e^\lambda\Theta)^{\tau}_0$.

Define $\gamma = \max(\{\alpha\}\cup \dom(r))$. Before defining $\eta$, observe that for all $\xi\in \dom(r)$, $r(\xi)<\ell^\xi e^\lambda \Theta$, so that by Lemma \ref{exponentialcancellation}, $e^\xi r(\xi)< e^\lambda \Theta=e^\gamma e^{\omega^\beta} \Theta$. Thus by normality of $e^\gamma$ and the fact that $e^{\omega^\beta} \Theta$ is a limit, for $\eta<e^{\omega^\beta}\Theta$ large enough we have that $e^{\gamma}\eta>e^{\xi} r (\xi)$ for all $\xi\in \dom(r)$. We claim that $U\cap (\eta,e^{\omega^\beta}\Theta]^{\tau}_\gamma\subset B^\tau_r(x) $.

To see this, assume that $x \not = y \in U\cap (\eta,e^{\omega^\beta}\Theta]^{\tau}_\gamma$ and let $\theta=\rho_\tau y$, so that $\theta \in (\eta,e^{\omega^\beta}\Theta]_\gamma$. Then, for $\xi\in \dom(r)$ we see that $\ell^\gamma\theta=\ell^{-\xi+\gamma}\ell^\xi\theta>\eta$ and hence $\ell^\xi\theta>e^{-\xi+\gamma} \eta$ by Lemma \ref{exponentialcancellation}. But $e^\xi e^{-\xi+\gamma} \eta=e^\gamma\eta>e^\xi r(\xi)$, so by normality of $e^\xi$, $e^{-\xi+\gamma} \eta> r(\xi)$ and thus $\ell^\xi\theta> r(\xi)$.

Meanwhile, for any $\xi<\lambda$, since $\theta<e^{\lambda}\Theta=e^\xi e^{-\xi+\lambda}\Theta$ by the way we chose $U$, it follows once again from Lemma \ref{exponentialcancellation} that $\ell^\xi\theta\leq e^{-\xi+\lambda} \Theta =\ell^\xi e^{\lambda}\Theta  $. Hence, for all $\xi\in \dom(r)$, $r(\xi)<\ell^\xi \theta <\ell^\xi e^{\lambda}\Theta$, that is, $\rho y=\theta\in(r(\xi),\ell^\xi e^\lambda\Theta]_{\xi}$ and thus $y\in B_r^\tau(x)$

We conclude that $U\cap (\eta, e^{\omega^\beta}\Theta]^{\tau}_\gamma\subset B_r^\tau(x)$, and thus $U\cap (\eta, e^{\omega^\beta}\Theta]^{\tau}_\gamma \subset V$, as needed. The modification for $\mathfrak{X}$ an ordinal with the left topology follows immediately as a special case.
\endproof

The following particular cases will be used later on:

\begin{cor}\label{nhbase3}
Let $\lambda$ be additively indecomposable. Then for any neighborhood $V$ of $e^\lambda\Theta$ in $[1, e^\lambda\Theta]_\lambda$ there exist ordinals $\eta < e^\lambda\Theta$ and $\zeta < \lambda$ such that $V$ contains the set $(\eta, e^\lambda\Theta]_{\zeta}$.
\end{cor}

Below, if $(X,\tau)$ is a scattered space and $(x_\xi)_{\xi<\mu}$ is a sequence  of points in $X$, we write $x_\xi\stackrel \lambda\to y$ if $x_\xi$ converges to $y$ in $\ic\lambda\tau$; that is, if for every $\lambda$-neighborhood $U$ of $y$, there is $\delta<\mu$ such that $x_\xi\in U$ whenever $\xi>\delta$.

\begin{lemma}\label{nhbase4}
Let $\lambda>0$, $\Theta$ be any ordinal and $\mu=\max(\cf \lambda,\omega)$. Then, there exists a sequence $(\theta_\xi)_{\xi<\mu}$ such that $\theta_\xi\stackrel \lambda\to e^\lambda(\Theta+1)$.
\end{lemma}

\proof
Write $\lambda=\alpha+\omega^\beta$. We will define $\theta_\xi$ for $\xi<\mu$ considering two cases. If $\beta=0$, then $\mu=\omega$. For $n<\omega$ define $\delta_n=\omega^\Theta\cdot n$, so that $(\delta_n)_{n<\omega}$ is  cofinal in $\omega^{\Theta+1}$. Then, define $\theta_n=e^\alpha \delta_n$, which is cofinal in $e^\lambda (\Theta+1)=e^{\alpha+1}(\Theta+1)=e^\alpha \omega^{\Theta+1}$ by normality of $e^\alpha$. 

If $\beta>0$, choose a sequence $(\gamma_\xi)_{\xi <\mu}$ which is cofinal in $\omega^\beta$. Then, set $\eta_\xi= e^{\gamma_\xi}(e^{\omega^\beta}\Theta+1)$, so that the sequence $(\eta_\xi)_{\xi <\mu}$ is cofinal in $e^{\omega^\beta}(\Theta+1)$ by Lemma \ref{recexp}.\ref{RecHypIV}. Finally, define $\theta_\xi=e^\alpha \eta_\xi$.

In either case, it is straightforward to check using Lemma \ref{nhbase2} that $\theta_\xi\stackrel\lambda\to e^\lambda(\Theta+1)$, as desired.
\endproof

Lemma \ref{nhbase2} and Corollary \ref{nhbase3} will be particularly useful for describing neighborhoods around points of high ranks in the completeness proof to follow.

\section{Completeness for generalized Icard spaces}\label{SecGL}

In this section, we will prove the strong completeness of $\gl$ with respect to its topological semantics. As we previously mentioned, we intend to construct $d$-maps from ordinal spaces onto $\omega$-bouquets. We start with a simple case:

\begin{lemma}\label{tm0}
For all nonzero additively indecomposable $\lambda$ and any $\omega$-bouquet generated by $(T,R)$ whose root $r$ has rank $\Theta$, there exists a surjective $d$-map $f:([0, e \Theta], \mathcal{I}_1) \to T$ such that $f(e\Theta)$ is the root of the tree.
\end{lemma}

\proof
By induction on $\Theta$. Denote by $r$ the root of the tree. If $\Theta = 0$, then the mapping $0 \mapsto r$ is clearly a $d$-map. Otherwise, suppose the result holds for all ordinals $\theta < \Theta$. We fix an enumeration $(v_i)_{i < \omega}$ of all daughters of $r$ with different properties, according as to whether $\Theta$ is a successor ordinal or a limit. If $\Theta$ is a limit ordinal, we let $(v_i)_{i < \omega}$ enumerate each daughter exactly once. If $\Theta$ is a successor ordinal, let $(v_i)_{i < \omega}$ be any enumeration where each daughter appears infinitely often. To each $v_i$ we associate its rank $\theta_i$ and its generated subbouquet $T_i$. Define $\alpha_0 = 0$ and, recursively, $\beta_i = \alpha_i + \omega^{\theta_i}$ and $\alpha_{i+1} = \beta_i + 1$.

If $\Theta=\Xi+1$ is a successor ordinal, then $\theta_i=\Xi$ for infinitely many $i$. If $\Theta$ is a limit, then the sequence $(\theta_i)_{i < \omega}$ is cofinal in $\Theta$. In either case it follows that $\beta_i\stackrel 1\to e\Theta$ and thus the sets $X_i := [\alpha_i, \beta_i]$ form a partition of $[0, e\Theta)$. 

By induction hypothesis, there are surjective $d$-maps $g_i: [0, e{\theta_i}] \to T_i$, and thus by Lemma \ref{disphomeomorphism} we obtain surjective $d$-maps $f_i\colon X_i\to T_i$ given by $f_i(\xi)=g_i(-\alpha_i+\xi)$. Let $f = \bigcup_{i < \omega}f_i \cup \{ (e\Theta, r) \}$. As the sets $X_i$ are clopen, $f$ is a $d$-map when restricted to $e\Theta$. We show it is also a $d$-map in $[0, e\Theta]$. 

Clearly, the function is pointwise discrete and surjective. Moreover, the family $\{[\alpha_i, e\Theta]\}_{i < \omega}$ forms a decreasing neighborhood base around $e\Theta$, and it is easy to check that $f[\alpha_n, e\Theta]=\{r\}\cup \bigcup_{i\geq n} T_i$, which is open, so $f$ is open. As for continuity, $f^{-1}T = [0, e\Theta]$, so we just need to check that the preimage of basic open sets around $r$ are open when $\Theta$ is a limit. But this is also true, as for any such set $A = \{r\} \cup \bigcup_{i \geq n}T_i$, we have $f^{-1}(A) = [\alpha_n, e\Theta]$ by construction.
\endproof

The $d$-maps constructed in Lemma \ref{tm0} can be further exploited. As the composition of $d$-maps is a $d$-map, it suffices to find $d$-maps from any space onto ordinals with the interval topology to obtain $d$-maps onto $\omega$-bouquets. As stated by Lemma \ref{LemmHyperRank}, this is sometimes possible by using logarithms $\ell^{\xi}\colon \Theta_{\xi + 1} \to \Theta_1$.

This is, however, impossible for limit ordinal Icard topologies. Instead, we could try to find different mappings. For example, the \emph{$\Lambda$-reductive maps} from \cite{polytopologies} function for limit ordinal Icard topologies as long as the index is countable. Theorem \ref{impossible $d$-maps} below will show that it is impossible to go further using this technique.

\begin{theorem}\label{impossible $d$-maps}
Let $\lambda$ be a nonzero ordinal, $\kappa$ be a limit ordinal, and $\Theta>0$ be any ordinal. If  $\max (cf(\kappa) , \aleph_0) \not= \max (cf(\lambda) , \aleph_0)$, then there exist no $d$-maps $f\colon (e^\kappa \Theta+1)_\kappa\to(e^\lambda \Theta+1)_\lambda$.
\end{theorem}

\proof Write $\lambda=\alpha+\omega^\beta$ and $\kappa=\gamma+\omega^\delta$, so that $\cf {\omega^\beta}=\cf \lambda$ and $\cf {\omega^\delta}=\cf \kappa$. It suffices to prove the theorem for $\Theta=1$ and $\lambda=\omega^\beta$, for if $f\colon (e^\kappa \Theta+1)_\kappa\to(e^\lambda \Theta+1)_\lambda$ were a $d$-map, we would have that $\ell^\lambda f (e^\kappa 1)=1$ (since $d$-maps are rank-preserving by Lemma \ref{propertiesofdmaps}.\ref{propertiesofdmapsOne}), and thus by Lemma \ref{LemmPropLog}.\ref{LemmPropLogIV}, there would exist $\sigma\in[\alpha,\lambda)$ such that $\ell^\sigma f(e^\kappa 1)=e^{\omega^\beta}1$. But then it is straightforward to check that $\tilde f\colon (e^\kappa 1+1)_\kappa\to(e^{\omega^\beta} 1+1)_{\omega^\beta}$, defined by $\tilde f(\xi)=\ell^\sigma f (\xi)$ if $\xi\in f^{-1}\ell^{-\sigma}[0,e^{\omega^\beta} 1]$, and $\tilde f(\xi)=0$ otherwise, would also be a $d$-map.

Let $\mu=\max (cf(\kappa) , \aleph_0)$ and $\nu=\max (cf(\lambda) , \aleph_0)$. By Corollary \ref{cofexp}, $\mu = \cf(e^\kappa 1)$ and $\nu = \cf(e^\lambda 1)$. Towards a contradiction, suppose that $\mu \neq \nu$ and $f\colon (e^\kappa 1+1)_\kappa\to(e^\lambda 1+1)_\lambda$ is a $d$-map. First assume that $\mu < \nu$. Use Lemma \ref{nhbase4} to find a sequence $(\theta_\xi)_{\xi<\mu}$ such that $\theta_\xi\stackrel\kappa\to e^\kappa 1$. By continuity, we should also have $f(\theta_\xi)\stackrel\lambda\to e^\lambda 1$; but, since $f$ preserves rank, $f(\theta_\xi)<e^\lambda 1$ for all $\xi$, and since $\mu < \nu = \cf {e^\lambda 1}$, we have that $(f(\theta_\xi))_{\xi<\mu}<\theta_\ast$ for some $\theta_\ast<e^\lambda 1$. But this means that $f(\theta_\xi)\not\in (\theta_\ast,e^\lambda 1]_1$ for all $\xi<\mu$, contradicting that $f(\theta_\xi)\stackrel\lambda\to e^\lambda 1$.

Now assume that $\mu > \nu$. Use Lemma \ref{nhbase4} once again to find a sequence $(\delta_\xi)_{\xi<\nu}$ such that $\delta_\xi\stackrel\lambda\to e^\lambda 1$. By continuity, $f^{-1}(\delta_\xi , e^\lambda 1]_0$ is $\lambda$-open for all $\xi<\nu$, hence by Corollary \ref{nhbase3}, there are sequences $(\eta_\xi)_{\xi<\nu}\subset e^\kappa 1$ and $(\gamma_\xi)_{\xi<\nu}\subset \kappa$ such that, for all $\xi<\nu$, $[0,e^\kappa 1]_0\cap (\eta_\xi,e^{\kappa}1]_{\gamma_\xi}\subset f^{-1}(\delta_\xi , e^\lambda 1]_0.$ Since $\nu < \mu = \cf{e^\kappa 1}$, there are $\eta_\ast<e^\kappa 1$ and $\gamma_\ast<\kappa$ bounding $(\eta_\xi)_{\xi<\nu}$ and $(\gamma_\xi)_{\xi<\nu}$, respectively. Then, $U_\ast=[0,e^\kappa 1]_0\cap (\eta_\ast,e^{\omega^\beta}1]_{\gamma_\ast}$ is a $\kappa$-neighborhood of $e^\kappa 1$, so that $f(U_\ast)$ contains a $\lambda$-neighborhood $V_\ast$ of $e^\lambda 1$. But, $V_\ast$ cannot contain any $\delta_\xi$, contradicting the fact that $\delta_\xi\stackrel\lambda\to e^\lambda 1$.

In either case, we conclude that there can be no such $d$-map.
\endproof

Our strategy will hence be to instead construct $d$-maps directly from Icard spaces $\ic\lambda{\mathfrak X}$ onto $\omega$-bouquets. The construction will be an adaptation of the one on Lemma \ref{tm0}, but the general case will be rather more involved. We will first assume that $\lambda$ is additively indecomposable and $\mathfrak X$ is an ordinal with the left topology for, as we shall see, the general case follows quickly from this particular one. As in Lemma \ref{tm0}, we will proceed by induction on the rank.

\subsection{The successor stage}\label{SuccStage} For this section, we fix an additively indecomposable ordinal $\lambda$ and an $\omega$-bouquet $(T,\sigma_R)$  with root $r$ of rank $\Theta + 1$. Let $(v_i)_{i<\omega}$ list all daughters of $r$ in such a way that each daughter is counted infinitely often, and let $T_i$ be the subbouquet generated by $v_i$. Denote the rank of $T_i$ by $\theta_i$. We introduce the following auxiliary notation:

\begin{defi} \label{alphas and betas}
Let $\iota < \lambda$. By Lemma \ref{ordinal arithmetic}.\ref{OrdArII}, there is a unique $k=k(\iota)<\omega$ such that $\iota$ can be written as $\alpha\cdot\omega + k$. For $\iota<\lambda$, define $\beta_{\iota} := e^{\iota + 1}(e^\lambda\Theta + 1 + e^\lambda \theta_{k(\iota )} )$. Then define $\alpha_\iota$ by cases, as $\alpha_0 := 0$, $\alpha_{\iota + 1} := \beta_\iota + 1$, and $\alpha_{\iota} := e^\iota(e^\lambda\Theta + 1)$ if $\iota$ is a limit.

Finally, define $X_\iota = [\alpha_\iota, \beta_\iota]$,  $Y_\iota = X_\iota\cap [0,1+e^\lambda\Theta]_{\iota+1}$, and $Z_\iota = X_\iota\setminus Y_\iota$.
\end{defi}


\begin{lemma}\label{its a partition}
The sets $X_\iota$,  $Y_\iota$, and $Z_\iota$ have the following properties:
\begin{enumerate}[(i)]

\item\label{LemmPartI} for all $\iota<\lambda$, each of the sets $X_\iota$, $Y_\iota$ and $Z_\iota$ are $\lambda$-clopen;

\item\label{LemmPartII} the collection $\{X_\iota \colon \iota < \lambda\}$ forms a partition of $e^\lambda(\Theta + 1)$;

\item\label{LemmPartZero} if $U$ is a $\lambda$-neighborhood of $e^\lambda(\Theta+1)$, then $Z_\iota\subset U$ for all $\iota$ large enough, and

\item\label{LemmPartIII} for all $\iota<\lambda$, $\ell^{\iota+1}(X_\iota)=[0,e^{\lambda}\Theta+1+e^{\lambda}\theta_{k(\iota )}],$ $\ell^{\iota+1}(Y_\iota)=[0,e^{\lambda}\Theta]$ and $\ell^{\iota+2}(Z_\iota)=[0,e^{\lambda}\theta_{k(\iota )}]$.
\end{enumerate}
\end{lemma}

\proof
Claim \eqref{LemmPartI} follows immediately from their definitions which use clopen intervals. For claim \eqref{LemmPartII}, note that from Lemma \ref{explimits}.\ref{explimits1} it follows that for limit $\iota$: $\lim_{\eta \to \iota} \alpha_\eta = \alpha_\iota$. Hence, the sets $X_\iota$ are disjoint and no gaps are left between them. Moreover, from Lemma \ref{explimits}.\ref{explimits2}, it also follows that $e^\lambda(\Theta+1) = \lim_{\iota \to \lambda}\alpha_\iota$.

For claim \eqref{LemmPartZero}, let $U$ be any $\lambda$-neighborhood of $e^\lambda(\Theta+1)$, so that by Lemma \ref{nhbase3} there are $\eta<e^\lambda(\Theta+1)$ and $\delta<\lambda$ such that $(\eta,e^\lambda(\Theta+1)]_\delta\subset U$. By Lemma \ref{recexp}.\ref{RecHypIV}, for $\gamma$ large enough we have that $ \eta< e^\gamma (e^\lambda \Theta+1)$. Since $\lambda$ is additively indecomposable, $\delta+\gamma<\lambda$, hence if $\iota>\delta+\gamma$ we have that for $\zeta \in Z _\iota$, $\ell^\gamma\ell^\delta\zeta= \ell^{\delta+\gamma}\zeta\geq \ell^{\iota+1}\zeta\geq e^{\lambda}\Theta+1$, so that by Lemma \ref{exponentialcancellation}, $\ell^\delta\zeta\geq e^\gamma (e^\lambda \Theta+1)>\eta$, that is, $\zeta\in (\eta,e^\lambda(\Theta+1)]_\delta\subset U$ and, since $\zeta$ was arbitrary, $Z _\iota\subset U$.

To show \eqref{LemmPartIII}, fix some $\iota < \lambda$. Since $\lambda$ is an additively indecomposable limit, $\iota+1+\lambda=\lambda$. Thus by Lemma \ref{LemmHyperRank}, $
\ell^{\iota + 1}\colon (\beta_\iota+1)_\lambda \to (e^{\lambda}\Theta+1+e^{\lambda}\theta_{k(\iota)}+1)_\lambda$ is a $d$-map. Also by Lemma \ref{LemmHyperRank}, $\ell^{\iota+1}$ is the rank function with respect to $\mathcal I_{\iota+1}$, and $X_\iota$ is $(\iota+1)$-open (in fact, $1$-open). Therefore, since $\beta_\iota\in X_\iota$ and $\ell^{\iota+1}\beta_\iota=e^{\lambda}\Theta+1+e^{\lambda}\theta_{k(\iota)}$, it follows from Lemma \ref{propertiesofdmapsTwo} that $[0,e^{\lambda}\Theta+1+e^{\lambda}\theta_{k(\iota)}]\subset \ell^{\iota + 1}(X_\iota)$. Since $Y_\iota=\{\xi\in X_\iota:\ell^{\iota+1}\xi\leq e^\lambda\Theta\}$, it follows that $\ell^{\iota+1}(Y_\iota)=[0,e^\lambda\Theta]$.

Finally, observe that if $\xi\in Z$, we can write $\ell^{\iota+1}\xi=e^\Lambda\Theta+1+\zeta$ for some $\zeta\leq 1+e^\lambda \theta_{k(\iota)}$, thus by Lemma \ref{LemmPropLog}.\ref{LemmPropLogII}, $\ell^{\iota+2}\xi=\ell (e^\Lambda\Theta+1+\zeta)=\ell(1+\zeta)\leq e^\lambda \theta_{k(\iota)}$, but $Z$ is $(\iota+2)$-open, so once again by Lemma \ref{propertiesofdmapsTwo}, $\ell^{\iota+2}(Z)=[0,e^{\lambda}\theta_{k(\iota )}].$
\endproof

The induction hypothesis for a successor stage will amount to assuming we are given fragments of a $d$-map, which we then need to complete. We make this notion precise now:

\begin{defi}\label{tm1}
We define a {\em pre-$d$-map for $(\Theta,\lambda)$ over $T$} to be a collection of $d$-maps $(g_i)_{i<\omega}$ such that, for every $i<\omega$, $g_i\colon (e^\lambda\theta_i+1)_\lambda\to T$ and $g_i$ has range $T_i$.
\end{defi}

\begin{lemma}\label{tm1}
Assume there exists a pre-$d$-map $(g_i)_{i<\omega}$ for $(\Theta+1,\lambda)$ over $T$. Then there is a surjective $d$-map $f\colon (e^\lambda(\Theta+1)+1)_\lambda \to T$.
\end{lemma}

\proof
We may assume without loss of generality that $\theta_0=\Theta$. As before, we have $k(\iota)$ denote the remainder of $\iota$ modulo $\omega$.

Before defining $f$, we set for $\iota<\lambda$, $f_\iota\colon X_\iota\to T_0\cup T_{k(\iota)}$ by $f_\iota(\xi) := g_0 \circ \ell^{\iota + 1}(\xi)$ if $\xi \in Y_\iota$, and $f_\iota(\xi) :=g_{k(\iota)} \circ \ell^{\iota + 2} (\xi)$ if $\xi \in Z_\iota.$ We claim that each $f_\iota: X_\iota \to T_0 \cup T_{{k(\iota)}}$ is a surjective $d$-map, with $f_\iota(Y_\iota)=T_0$ and $f_\iota(Z_\iota)=T_{k(\iota)}$:

\noindent{\sc Surjectivity.} That $f_\iota(Y_\iota)=T_0$ and $f_\iota(Z_\iota)=T_{k(\iota)}$ follows easily from \ref{its a partition}.\ref{LemmPartIII} and the assumption that each $g_i$ is onto $T_i$, and hence $f_\iota$ is surjective.

\noindent{\sc Openness.} Let $U \subset X_\iota$ be open. Define $U_1=U\cap Y_\iota$, and $U_2=U\setminus U_1\subset Z_\iota$, both of which are open. Because $g_0$, $g_\iota$, $\ell^{\iota+1}$ and $\ell^{\iota+2}$ are all $d$-maps, $f_\iota(U)$ is the union of the open sets $f_\iota(U_1)$ and $f_\iota(U_2)$, so $f_\iota(U)$ is open; since $U$ was arbitrary, $f_\iota$ is an open function.

\noindent{\sc Continuity.} This follows by a similar argument: if $V\subset T_0\cup T_{k(\iota)}$ is open, then $f^{-1}_\iota(V)=(f^{-1}_\iota (V)\cap Y_\iota) \cup (f^{-1}_\iota (V)\cap Z_\iota)$, which is a union of two open subsets of $X_\iota$.

\noindent{\sc Pointwise discreteness.} If $t \in T_0\cup T_{k(\iota)}$, then once again $f_{\iota}^{-1}(t)=(f^{-1}_\iota (V)\cap Y_\iota) \cup (f^{-1}_\iota (V)\cap Z_\iota)$. Since $g_0 \circ \ell^{\iota + 1}$ is a $d$-map, $A=f^{-1}_\iota (V)\cap Y_\iota$ is a discrete subset of $Y_\iota$, and similarly $B=f^{-1}_\iota (V)\cap Z_\iota$ is a discrete subset of $Z_\iota$. But both $Y_\iota$ and $Z_\iota$ are clopen, so $A\cup B$ is also discrete.

This completes the proof of the claim. We define $f=\bigcup_{\iota \in \lambda}f_\iota \cup \{(e^{\lambda}(\Theta+1),r)\}$. It can be quickly verified that $f$ is surjective. We show it is also a $d$-map:

\noindent{\sc Openness.} If $U$ is open in $[0,e^{\lambda}(\Theta+1))$, then $f(U)=f(\bigcup_{\iota \in \lambda}(U \cap X_\iota)=\bigcup_{\iota \in \lambda} f_\iota(U \cap X_\iota)$ is open since each $f_\iota$ is open. Otherwise, if $U$ is a neighborhood of $e^{\lambda}(\Theta+1)$, we claim that $f(U)=T$. Indeed, $r=f(e^\lambda(\Theta+1))\in f(U)$, and if $t\not=r$, we have that $t\in T_i$ for some $i$. By Lemma \ref{its a partition}.\ref{LemmPartZero} we have that, for some $\iota_\ast<\lambda$, $Z_\iota\subset U$ whenever $\iota>\iota_\ast$; in particular, since every subbouquet is counted infinitely often and $\lambda$ is a limit, we can take $\iota>\iota_\ast$ so that $k(\iota)=i$, and since $T_i=f(Z_\iota)$, we have that $t\in f(Z_\iota)\subset f(U)$. It follows that $f$ is open.

\noindent{\sc Continuity.} Suppose $V\in \sigma_R$ and $\xi\in f^{-1}(V)$. If $f(\xi)=r$, then since the only neighborhood of $r$ is all of $T$, it follows that $f^{-1}(V)=e^{\lambda}(\Theta+1)+1$. Otherwise, by Lemma \ref{its a partition}.\ref{LemmPartII}, $\xi\in X_\iota$ for some unique $\iota$. Then, $\xi\in f^{-1}(V)\cap X_\iota\subset f^{-1}(V);$ but, $f^{-1}(V)\cap X_\iota= f^{-1}_\iota(V\cap T_{k(\iota)})$, which is open since $f_\iota$ is continuous. Since $\xi\in f^{-1}(V)$ was arbitrary, we conclude that $f^{-1}(V)$ is open.

\noindent{\sc Pointwise discreteness.} Each $f_\iota^{-1}(t)$ is discrete and each $X_\iota$ is $(\iota+1)$-clopen and thus ${\lambda}$-clopen, so $f^{-1}(t)$ is discrete for each $t\not=r$, while $f^{-1}(r)$ is a singleton and thus discrete.

Therefore, $f$ is a $d$-map, as desired.
\endproof

Note that the above may already be used to give an inductive construction of $d$-maps onto any $\omega$-bouquet of finite rank, but we have yet to deal with limit $\Theta$.

\subsection{The limit stage} For the following, we fix an additively indecomposable ordinal $\lambda$ and an $\omega$-bouquet $(T, \sigma_R)$ with root $r$ of countable limit rank $\Theta$. We also fix a sequence $(v_i)_{i < \omega}$ enumerating each daughter of $r$ exactly once and denote by $T_i$ and $\theta_i$, respectively, the generated subbouquet and rank of $v_i$.  

A {\em dominating subsequence of $(\theta_i)_{i<\omega}$} is a subsequence $(\theta_{m_i})_{i<\omega}$ of $(\theta_i)_{i<\omega}$ such that $(\theta_{m_i})_{i<\omega}$ is strictly increasing and for all $i<\omega$, $\theta_i<\theta_{m_i}$. Observe that if this is the case, then $\lim_{i\rightarrow \omega}\theta_{m_i}\rightarrow \Theta$. It is evident that dominating sequences for $(\theta_i)_{i<\omega}$ do exist. We fix one such sequence $(\theta_{m_i})_{i<\omega}$.

\begin{lemma}\label{LemmLimitInductiveBasic}
If there exists a pre-$d$-map $(g_i)_{i<\omega}$ for $(\Theta,\lambda)$ over $T$, then there is a continuous surjective map $h\colon (e^\lambda\Theta+1)_\lambda\to T$ such that $h^{-1}(r)=\{e^\lambda\Theta\}$, and $h|_{e^\lambda\Theta}$ is a $d$-map.
\end{lemma}

\proof
%
Let $h\colon e^\lambda \Theta+1\to T$ be defined by $h(e^\lambda\Theta)=r$, $h(\xi)=g_{m_0}(\xi)$ if $\xi\leq e^\lambda\theta_{m_0}$, and $h(\xi)=g_{m_i}(\xi)$ if $e^\lambda\theta_{m_{i-1}}<\xi\leq e^\lambda\theta(m_{i})$ with $0 < i$. It is straightforward to check that $h$ has the desired properties.
\endproof

The map $h$ just defined may well not be a $d$-map, as it might fail to be open. In the following, we will modify $h$ so as to make it open near $e^\lambda\Theta$. This will require constructing some auxiliary sets, in a fashion similar to Definition \ref{alphas and betas}.

\begin{defi}
Let $k(\iota)$ denote the remainder of $\iota$ modulo $\omega$ as before. For $j<\omega$ and $\iota<\lambda$, we set $\gamma_{j \iota}= e^\lambda \theta_{m_{j+1}} + e^\iota ( e^\lambda \theta_{m_{j}}+1)$ and $\delta_{j \iota}= \gamma_{j\iota}+e^\iota \big (e^\lambda\theta_{m_{j}}+ 1 + e^\lambda \theta_j \big )$.

We also define $W_{j\iota}=(\gamma_{j\iota},\delta_{j\iota}]_0 \cap (e^\lambda\theta_{m_j},e^\lambda\theta_{m_{j+1}}+ 1 + e^\lambda \theta_j]_\iota$ and $W_j=\bigcup_{\iota<\lambda}W_{j\iota}$.

\end{defi}

\begin{lemma}\label{LemmW}
The sets $W_{j\iota}$ and $W_j$ have the following properties:
\begin{enumerate}[(i)]

\item\label{LemmWI} $\ell^{\iota+1} ( W_{j\iota})= e^\lambda\theta_j+1$ for all $j<\omega$ and $\iota<\lambda$,

\item\label{LemmWII} $W_j\subset (e^\lambda \theta_{m_{j+1}},e^\lambda \theta_{m_{j+2}})$ and $W_{i\kappa}\cap W_{j\iota}=\varnothing$ if $i\not=j$ or $\iota\not=\kappa$,

\item\label{LemmWIII} the set $W_j$ is closed for all $j$, and

\item\label{LemmWIV} if $U$ is any $\lambda$-neighborhood of $e^\lambda\Theta$, there exist $n<\omega$ and $\iota<\lambda$ such that $W_{j\iota}\subset U$ for all $j>n$.

\end{enumerate}

\end{lemma}

\proof
Claim \eqref{LemmWI} is similar to Lemma \ref{its a partition}.\ref{LemmPartIII}.
We prove claim \eqref{LemmWII}. Observe that for any $j<\omega$, $\gamma_{j\iota}<\delta_{j\iota}<\gamma_{j\kappa}$ whenever $\iota<\kappa$, from which it follows that $W_{j\kappa}\cap W_{j\iota}=\varnothing$ if $\iota\not=\kappa$. Moreover, if we define $\nu_j=e^{\lambda}\theta_{m_{j+1}}+e^{\lambda}(\theta_{m_{j}}+1)< e^{\lambda}\theta_{m_{j+2}}$, by Lemma \ref{explimits}.\ref{explimits1} we have that $\gamma_{j\iota},\delta_{j\iota} \stackrel 1\to \nu_j $ as $\iota\to \lambda$ , which by the definition of $\gamma_{j 0}$ implies that $W_j\subset (e^\lambda \theta_{m_{j+1}},\nu_j) \subset (e^\lambda \theta_{m_{j+1}},e^\lambda \theta_{m_{j+2}})$. It follows that $W_{i\kappa}\cap W_{j\iota}=\varnothing$ if $i\not=j$.

For claim \eqref{LemmWIII}, since each $W_{j\iota}$ is closed, it also follows that $\overline W_j\setminus W_j\subset \{\nu_j\}$. However, $[0,\nu_j]_0\cap (e^{\lambda} \theta_{m_{j+1}},\infty)_1$ is a $\lambda$-neighborhood of $\nu_j$ which does not intersect any of the intervals $(\gamma_{j\iota},\delta_{j\iota}]$. We conclude that $\overline W_j\setminus W_j$, that is, $W_j$ is closed.

Finally, for claim \eqref{LemmWIV}, let $U$ be any neighborhood of $e^\lambda\Theta$, so that by Lemma \ref{nhbase3} and the fact that $\theta_{m_i}\to \Theta$, there are $n<\omega$ and $\iota<\lambda$ such that $[0,e^\lambda\Theta]_0 \cap (e^\lambda\theta_{m_n},e^\lambda\Theta]_\iota\subset U$. Then, it is immediate that if $j>n$, $W_{j\iota}\subset U$.
\endproof

\begin{lemma}\label{LemmLimitInductiveStrong}
Assume there exists a pre-$d$-map $(g_i)_{i<\omega}$ for $(\Theta,\lambda)$ over $T$. Then there is a surjective $d$-map $f\colon ( e^\lambda\Theta+1 )_\lambda \to T$.
\end{lemma}

\proof
We use Lemma \ref{LemmLimitInductiveBasic} to construct a continuous map $h\colon (e^\lambda\Theta+1)_\lambda\to T$ such that $h^{-1}(r)=\{e^\lambda\Theta\}$, and $h|_{e^\lambda\Theta}$ is a $d$-map. Then, we define $f$ by $f(\xi)=h(\xi)$ if $\xi\not\in\bigcup_{j<\omega}W_j$ and $f(\xi)=g_j\ell^{\iota+1}\xi$ if $\xi\in W_{j\iota}$. That $f$ is surjective and pointwise discrete is clear from its construction. We verify openness and continuity.

\noindent{\sc Openness.} Choose $\xi < e^\lambda\Theta$. If $\xi\in W_{j\iota}$ for some $j,\iota$, we use the fact that $W_{j\iota}$ is open and that $g_j$, $\ell^{\iota+1}$ are $d$-maps to see that $f$ is open near $\xi$. Otherwise, $\xi\in (e^\lambda\theta_{m_j},e^\lambda\theta_{m_{j+1}}]$ for some $j$, and $(e^\lambda\theta_{m_j},e^\lambda\theta_{m_{j+1}}]\setminus \bigcup_{i<\omega} W_i= (e^\lambda\theta_{m_j},e^\lambda\theta_{m_{j+1}}]\setminus W_j$. But $W_j$ is $\lambda$-closed by Lemma \ref{LemmW}.\ref{LemmWIII}, whereby $(e^\lambda\theta_{m_j},e^\lambda\theta_{m_{j+1}}]\setminus W_j$ is $\lambda$-open. Since $h$ is open near $\xi$ then so is $f$. Hence, $f$ is open in $[0, e^\lambda\Theta)$.

Now let $V$ be an open neighborhood of $e^\lambda\Theta$. By Lemma \ref{LemmW}.\ref{LemmWIV}, there are $\iota<\lambda$ and $n<\omega$ such that $W_{j\iota}\subset U$ for all $j>n$, and it follows by the surjectivity of $g_j$ that $\{r\}\cup\bigcup_{j>n} T_j\subset f(V)$. Since $V$ was arbitrary, we conclude that $f$ is open near $e^\lambda\Theta$.

\noindent{\sc Continuity.} By a similar argument to the one above, $f$ is continuous in $[0, e^\lambda\Theta)$, and so it remains to show that it is continuous at $e^\lambda\Theta$. Take a neighborhood $U$ of $r$, so that $U=\{r\}\cup \bigcup_{i>n}T_i$, and choose $k>n$ large enough so that $m_k>n$. It is straightforward to check that $(e^\lambda \theta_{m_k+1},e^\lambda\Theta]\subset f^{-1}(U)$, so the latter is open near $e^\lambda\Theta$. 
\endproof

\subsection{Main result} First, we join Lemmas \ref{tm1} and \ref{LemmLimitInductiveStrong} to build $d$-maps from Icard spaces onto $\omega$-bouquets.

\begin{theorem}\label{TheoMapsBouq}
Let $(T,\sigma)$ be any $\omega$-bouquet with root $r$ of (countable) rank $\Theta$, and let $\lambda>0$ be any ordinal. Then, there is a surjective $d$-map $f\colon (e^\lambda\Theta+1)_\lambda\to T$.
\end{theorem}

\proof
If $\lambda=\alpha+\omega^\beta$ and $f\colon (e^{\omega^\beta}\Theta+1)_{\omega^\beta}\to T$ is a surjective $d$-map, then $f\circ \ell^\alpha \colon (e^{\lambda}\Theta+1)_{\lambda}\to T$ is also a surjective $d$-map. Hence, it suffices to consider additively indecomposable $\lambda$. If $\lambda=1$ the claim becomes Lemma \ref{tm0}, so we assume that $\lambda$ is a limit and proceed by induction on $\Theta$.

If $\Theta=0$, then $T=\{r\}$, and $f\colon e^\lambda 0+1\to T$ given by $f(0)=r$ trivially has all desired properties.

Otherwise, we may enumerate all immediate subbouquets of $T$ by $(T_i)_{i<\omega}$ and denote their ranks by $\theta_i$. By induction on their rank, we find surjective $d$-maps $g_i\colon (e^\lambda\theta_i+1)_\lambda\to T_i$, so that $(g_i)_{i<\omega}$ is a pre-$d$-map. Thus, by Lemma \ref{tm1} if $\Theta$ is a successor, or Lemma \ref{LemmLimitInductiveStrong} if $\Theta$ is a limit, there is a surjective $d$-map $f\colon (e^\lambda\Theta+1)_\lambda\to T$, as needed.
\endproof

\begin{cor}\label{tm1c}
Let $\lambda$ be a nonzero ordinal and $\Theta$ be countable. For each scattered space $\mathfrak X$ of rank $e^\lambda \Theta + 1$ and any $\omega$-bouquet $T$ of rank $\Theta$ there exists a surjective $d$-map $f:\ic\lambda{\mathfrak X} \to T$.
\end{cor}

\proof
From Corollary \ref{lambdadmaps} and Theorem \ref{TheoMapsBouq}.
\endproof

Note that it follows from Lemma \ref{propertiesofdmaps}.\ref{propertiesofdmapsOne} that if $\mathfrak X$ contains only one point of maximal rank, it is also the preimage of the root. Meanwhile, as a consequence of Theorem \ref{TheoMapsBouq}, we obtain the following:

\begin{theorem}[Strong completeness]\label{icardcompleteness}
Let $\lambda$ be a nonzero ordinal and $\mathfrak X$ be a scattered space of rank at least $e^\lambda \omega +1$. Then, $\gl$ is strongly complete with respect to $\ic\lambda{\mathfrak X}$.
\end{theorem}

\proof Suppose $\Gamma$ is a consistent set of formulae. By Theorem \ref{TheoCompBouq}, $\Gamma$ is satisfiable on some $\omega$-bouquet $T$ of rank $\Theta \leq \omega$, whence by Lemma \ref{propertiesofdmaps}.\ref{valuation transfer} and Theorem \ref{TheoMapsBouq}, $\Gamma$ is satisfiable on $(e^\lambda\Theta + 1)_\lambda$, a $\lambda$-open set of $\ic\lambda{\mathfrak X}$.
\endproof

As an immediate consequence of Theorem \ref{icardcompleteness}, we obtain the following result that is of particular importance when $\lambda$ is uncountable:

\begin{cor}\label{corints}
$\gl$ is strongly complete with respect to an ordinal $(\Theta, 
\ico \lambda)$ whenever $e^\lambda \omega <\Theta$.
\end{cor}

The instance of Corollary \ref{corints} when $\lambda = 1$ is a strenghtening of the Abashidze-Blass theorem. Another remarkable consequence of Theorem \ref{icardcompleteness} is what results when applying it to the club topology $\tau_c$, (see Example \ref{examples of scattered spaces}) as it has been shown (see \cite{blass1990}) that it is consistent with $\zfc$ + ``there exists a Mahlo cardinal'' that $\gl$ be incomplete with respect to $\tau_c$ for any ordinal. However, using generalized Icard topologies we obtain the following:

\begin{cor}\label{corclubs}
$\gl$ is strongly complete with respect to an ordinal $(\Theta, \ic \lambda {\tau_c})$ whenever $\aleph_{e^\lambda \omega + 1} < \Theta$.
\end{cor}

\section{Concluding remarks}

We have seen that $\gl$ is strongly complete with respect to the Icard topologies of any scattered space of sufficiently large rank. This is a rather remarkable property as, more frequently than not, $\gl$ is not complete with respect to the original space. For example, the space $\ico 0$ cannot satisfy the formula $\Diamond (p \land \Box \bot) \land \Diamond (\neg p \land \Box \bot),$ which is consistent with $\gl$, but any other $\ico \lambda$ can. An analogous situation occurs with the club topology, as seen in Corollary \ref{corclubs}, and with other topologies such that the consistency strength of the completeness of $\gl$ with respect to them is not even known, such as the so-called Mahlo topology (see \cite{bekgab14}) or the topology induced by the measurable filter (see \cite{blass1990}), although in those cases, the existence of points of sufficiently large rank also requires assumptions well beyond $\zfc$. 


Our construction relies heavily on the fact that the set $\mathbb{P}$ of propositional variables of $\gl$ is countable. Theorem \ref{icardcompleteness} may fail if this is not the case. As a simple example, assume $\vert \mathbb{P} \vert = (2^{\aleph_0})^+$ and $\Gamma = \{\Diamond p \colon p \in \mathbb{P}\} \cup \{\Box\neg(p \land q) \colon p,q \in \mathbb{P}\}$. On any countable space, there are always two variables that receive the same valuation, whereby $\Gamma$ cannot be satisfied. Nonetheless, it is easy to find generalized Icard spaces of higher cardinality that satisfy $\Gamma$. This gives rise to the following question:

\begin{question}
Assume $\gl$ is endowed with a set of propositional variables of cardinality $\kappa \geq \aleph_1$. Is there a natural topological space $\mathfrak X$ with respect to which $\gl$ is strongly complete?
\end{question}

As mentioned on the introduction, Icard topologies find a natural application in the construction of models of the polymodal logics $\glp_\Lambda$. Although the present work provides a significant advance towards extending known results, the completeness of $\glp_\Lambda$ for uncountable $\Lambda$ remains, as to now, unsettled.


\begin{thebibliography}{10}

\bibitem{abashidze1985}
M.~Abashidze.
\newblock Ordinal completeness of the {G}\"odel-{L}\"ob modal system.
\newblock {\em Intensional Logics and the Logical Structure of Theories}, pages
  49--73, 1985.
\newblock in Russian.

\bibitem{topocompletenessofglp}
L.~D. Beklemishev and D.~Gabelaia.
\newblock {Topological completeness of the provability logic GLP}.
\newblock {\em Annals of Pure and Applied Logic}, 164(12):1201--1223, 2013.

\bibitem{bekgab14}
L.~D. Beklemishev and D.~Gabelaia.
\newblock {Topological Interpretations of Provability Logic}.
\newblock {\em Leo Esakia on Duality in Modal and Intuitionistic Logics}, pages
  259--290, 2014.

\bibitem{bmm}
G.~Bezhanishvili, R.~Mines, and P.~Morandi.
\newblock {Scattered, Hausdorff-reducible, and hereditarily irresolvable
  spaces}.
\newblock {\em Topology and its applications}, 132(3):291--306, 2003.

\bibitem{bezhanishvilimorandi}
G.~Bezhanishvili and P.~J. Morandi.
\newblock {Scattered and hereditarily irresolvable spaces in modal logic}.
\newblock {\em Archive for Mathematical Logic}, 49:343--365, 2010.

\bibitem{blass1990}
A.~Blass.
\newblock Infinitary combinatorics and modal logic.
\newblock {\em Journal of Symbolic Logic}, 55(2):761--778, 1990.

\bibitem{esakia1981}
L.~Esakia.
\newblock Diagonal constructions, l\"ob's formula and cantor's scattered space
  (in russian).
\newblock {\em Studies in logic and semantics}, 132(3):128--143, 1981.

\bibitem{polytopologies}
D.~Fern\'andez-Duque.
\newblock The polytopologies of transfinite provability logic.
\newblock {\em Archive for Mathematical Logic}, 53(3-4):385--431, 2014.

\bibitem{modelsofglp}
D.~Fern\'andez-Duque and J.~J. Joosten.
\newblock Models of transfinite provability logic.
\newblock {\em Journal of Symbolic Logic}, 78(2):543--561, 2011.

\bibitem{hyperations}
D.~Fern\'andez-Duque and J.~J. Joosten.
\newblock {Hyperations, Veblen progressions, and transfinite iteration of
  ordinal functions}.
\newblock {\em Annals of Pure and Applied Logic}, 164(7-8):785--801, 2013.

\bibitem{godel31}
K.~G\"odel.
\newblock \"{U}ber formal unentscheidbare {S}\"atze der {P}rincipia
  {M}athematica und verwandter {S}ysteme {I}.
\newblock {\em Monatshefte f\"ur Mathematik Physik}, 38:173--198, 1931.

\bibitem{icardglp}
T.~F. Icard.
\newblock A topological study of the closed fragment of \glp.
\newblock {\em Journal of Logic and Computation}, 21(4):683--696, 2011.

\bibitem{jech}
T.~Jech.
\newblock {\em {Set Theory}}.
\newblock Springer monographs in Mathematics, 2006.

\bibitem{segerberg1971}
K.~Segerberg.
\newblock An essay in classical modal logic.
\newblock {\em Filosofiska F\'oreningen och Filosofiska Institutionen vid
  Uppsala Universitet}, 1971.

\bibitem{shehtman99}
V.~Shehtman.
\newblock On strong neighbourhood completeness of modal and intermediate
  propositinoal logics (part ii).
\newblock {\em JFAK. Essays Dedicated to Johan van Benthem on the Occasion of
  his 50th Birthday}, 1999.

\bibitem{jvblogicsspace}
J.~van Benthem and G.~Bezhanishvili.
\newblock {\em Modal logics of space}.
\newblock Institute for Logic, Language and Computation, 2006.

\end{thebibliography}
\end{document}